\documentclass{amsart}

\usepackage{graphicx}
\usepackage{amsmath} 
\usepackage{amssymb}  
\usepackage{graphics}
\usepackage{epsfig}
\usepackage[ruled,vlined]{algorithm2e}
\usepackage{subfig}
\usepackage{amsfonts}
\usepackage{mathrsfs}
\usepackage{bbm}
\usepackage{color}
\usepackage{pgfpages}
\usepackage{pgf, tikz}
\usepackage{tkz-berge}
\usetikzlibrary{arrows, automata, shapes, snakes,backgrounds}
\usepackage{comment}

\usepackage{tikz}
\usetikzlibrary{arrows.meta}
\usetikzlibrary{trees}

\newcommand{\beq}{\begin{equation}}
\newcommand{\eeq}{\end{equation}}

\newcommand{\E}{{\mathcal E}}

\newcommand{\R}{\mathbb R}

\def\e{{\rm u}}

\def\eps{\varepsilon}

\newcommand{\Id}{{\rm I}}

\newcommand{\epstildem}{\widetilde\eps_m}

\newcommand{\Ginv}{M^\#}
\newcommand{\Ainv}{A_\eps^\#}

\newcommand{\G}{\mathcal{G}}

\newcommand{\Rm}{\varrho_m}

\newcommand\sym{\mathrm{Sym}}
\newcommand\skw{\mathrm{Skew}}
\newcommand\Pact {P^+}

\newcommand{\vp}{\mathrm{v}}
\newcommand{\up}{\mathrm{u}}
\newcommand{\Gf}{R}

\newtheorem{lemma}{Lemma}[section]
\newtheorem{theorem}{Theorem}[section]

\newtheorem{defi}{Definition}[section]

\newtheorem{assumption}{Assumption}[section]

\begin{document}

\date{}

\title[Changing ranking in a graph]{Changing the ranking in eigenvector centrality of a weighted graph by small perturbations}

\author{Michele Benzi}

\address{Scuola Normale Superiore, Piazza dei Cavalieri 7, Pisa, Italy.}

\email{michele.benzi@sns.it}

\thanks{Michele Benzi's work was supported in part by MUR (Ministero dell' Universit\`a e della Ricerca) through the
PNRR MUR Project PE0000013-FAIR and through the PRIN Project 20227PCCKZ (``Low Rank Structures and Numerical
Methods in Matrix and Tensor Computations and their Applications").}

\author{Nicola Guglielmi} 

\address{Gran Sasso Science Institute (GSSI), via Crispi 7,
I-67010    L' Aquila,  Italy.}

\email{nicola.guglielmi@gssi.it}

\thanks{
Nicola Guglielmi acknowledges support from the Italian 
MUR within the 
PRIN 2022 Project ``Advanced Numerical Methods for Time Dependent Parametric 
Partial Differential Equations with Applications''.  
He also thanks the support from the MIUR program for Departments of Excellence
(project title of GSSI: Pattern analysis and engineering). \\
\indent
Both authors acknowledge the MUR Pro3 joint project 
entitled ``STANDS: Stability of Neural Dynamical Systems". 
The authors are also affiliated to the INdAM-GNCS (Gruppo Nazionale di Calcolo Scientifico).}



\begin{abstract}
In this article, we consider eigenvector centrality for the nodes of a graph and study the
robustness (and stability) of this popular centrality measure.
For a given weighted graph $\G$ (both directed and undirected), we consider the associated 
weighted adjacency matrix $A$, which by definition is a non-negative matrix.  
The eigenvector centralities of the nodes of $\G$ are the entries
of  the  Perron eigenvector of $A$, which is the (positive) eigenvector associated with the eigenvalue with 
largest modulus. They provide a ranking of the nodes according to the corresponding centralities. 
An indicator of the robustness of eigenvector centrality consists in looking for a nearby 
perturbed graph $\widetilde{\G}$, with the same structure as $\G$ (i.e., with the same 
vertices and edges), but with a weighted adjacency matrix $\widetilde A$ such that the 
highest $m$ entries ($m \ge 2$) of the Perron eigenvector of $\widetilde A$ coalesce, 
making the ranking at the highest level ambiguous. 
To compute a solution to this matrix nearness problem, a nested iterative algorithm is 
proposed that  makes use of a constrained gradient system of matrix differential 
equations 
in the inner iteration and a one-dimensional 
optimization of the perturbation size in the outer iteration.

The proposed algorithm produces the {\em optimal} perturbation (i.e., the one with 
smallest Frobenius norm) of the $A$ which causes the looked-for coalescence,  which is a 
measure of the sensitivity of the graph.
As a by-product of our methodology we derive an explicit formula for the gradient of
a suitable functional corresponding to the minimization problem, leading to a characterization
of the stationary points of the associated gradient system. This in turn gives us important properties 
of the solutions, revealing an underlying low rank structure. Our numerical experiments indicate
that the proposed strategy outperforms more standard approaches based on algorithms
for constrained optimization.
The methodology is formulated in terms of graphs but applies to any nonnegative 
matrix, with potential applications in fields like population models, consensus dynamics, 
economics, etc.
\end{abstract}



\pagestyle{myheadings}
\thispagestyle{plain}
\markboth{}{Robustness of eigenvector centrality measures}

\maketitle

\section{Introduction}

In many fields, such as biology, social sciences, statistics, medicine and several others, ranking 
entities in a large data set is an important task. This leads to the concept of centrality
measure for a network or a graph, which is the subject of this paper. High-rank (centrality)
of a node in a graph indicates the importance of the node. Several centrality 
measures have been proposed in the literature; for a broad class of walk-based centrality 
measures see, e.g, \cite{BeKl15,BeBo20,New10}); here we consider eigencentrality, i.e., the 
one based on the Perron eigenvector, which was first introduced by Landau \cite{Landau1895}
and rediscovered by several authors since.
For a strongly connected graph $\G$ with associated (possibly weighted) adjacency matrix $A$, this consists in ranking
the nodes according to the entries of the Perron eigenvector of $A$, which is the positive eigenvector
associated to the eigenvalue $\lambda>0$ of maximum modulus. 
As is well known, one the most important ranking schemes,
 the PageRank algorithm \cite{BrinPage1998,LM96}, is based on a variant of this idea.

Given a directed or undirected weighted graph, we are interested 
in understanding how small perturbations may affect the ranking
and which edges are more sensitive in order to produce certain changes
in the ranking. 
We can distinguish two main motivations.
\begin{itemize}
\item[(i) ] Perturbing the weights as a robustness tool, or sensitivity analysis.
Motivation: studying the reliability of the ranking in
presence of uncertainties or variations in time of the weights, 
which in real networks are never known exactly.
\item[(ii) ] Assessing the vulnerability of a network in the presence 
of deliberate attempts, for example by malicious, disgruntled individuals
or by other interested agents, to manipulate the rankings by slightly perturbing the weights.
E.g, in PageRank, one might want to manipulate the network
in order to affect the rank of a product or of a political candidate.
\end{itemize}

A robustness indicator  is useful in assessing the reliability of the measure. A very cheap and 
intuitively reasonable indicator of robustness consists of measuring the distance between the
first and the $m$-th largest entries of the (suitably normalized) ranking eigenvector. However 
such a measure may be misleading, which means that even when this distance is small, this does
not imply that the ranking is not robust, because - due to the structure of the graph - even a large
modification of the weights of the graph might not be sufficient to annihilate the considered distance.
For this reason we pose the problem in different terms and look for a closest graph $\widetilde{\G}$ to $\G$, which
has in common with $\G$ the same sets of vertices and edges, but different weights, such that the largest entries 
of the Perron eigenvector coalesce. This problem translates into the search of a perturbed weighted adjacency matrix 
$\widetilde A$ to $A$ which preserves the sparsity pattern of $A$ and certain of its properties, and is characterized by 
coalescence of the first $m$ entries ($m \ge 2$) of the leading eigenvector.

In this way we define a structured robustness measure, whose computation will be discussed in this paper, together
with the computation of the closest altered graph. 
It may be computationally expensive to be used as a general-purpose robustness indicator in comparison with simpler heuristic
criteria, which are available without further computations, such as the above mentioned difference between the largest and $m$-th largest entry
in the Perron eigenvector of $A$. But it has the advantage of a higher reliability. 
Moreover, cheap measures are not able to provide the nearby graph with ambiguous or changed ranking, which is a very
important byproduct of the methodology we propose.
The quantity considered in this paper is related to other structured robustness measures that arise in matrix 
and control theory, such as  stability radii, 
passivity distances,  distance to singularity, etc. 
\cite{HinrichsenPritchard86,Higham88,Higham02,BennerVoigt14,GKL15,GLM17,KressnerVoigt15}.

A main objective of this paper is to show how the robustness measure can be computed through the
explicit computation of a nearby graph with ambiguous ranking. 
The proposed algorithm is an iterative algorithm, where in each step the leading eigenvector of 
the weighted adjacency matrix of a graph with perturbed weights is computed. 
For large sparse graphs  with $n$ vertices and $O(n)$ edges, 
these computations can typically be done with a more favorable complexity in the number of vertices 
with respect to a general purpose eigensolver. 
A feature of this algorithm in common with recent algorithms for eigenvalue optimization 
(such as the ones described in \cite{GKL15,GLM17,AEGL19}) is a two-level methodology for matrix 
nearness problems, where in an 
inner iteration a gradient flow drives perturbations of a fixed norm to the original matrix  into 
a (local) minimum of a nonnegative functional that depends on eigenvalues and/or eigenvectors, while in an 
outer iteration the perturbation size is optimized such that the functional becomes zero. 

The paper is organized as follows. In Section~\ref{sec:related} we briefly mention recent related work
by other authors working on similar centrality manipulation problems.
In Section~\ref{sec:problem} we recall basic concepts of eigenvector centrality measures. 
In Subsection~\ref{sec:robustness measure} we introduce the radius of robustness ($\Rm$) for 
eigenvector centrality as the distance to the closest graph with coalescent
highest $m$ entries in the associated Perron eigenvector.
In Subsection~\ref{sec:two-level method} we formulate a methodology to compute the $\Rm$ and the
nearby ambiguous graph. 
This requires the solution of a structured matrix nearness problem. 
For this we propose a two-level iterative method, which is based on a gradient system of 
matrix differential equations in an inner iteration and a one-dimensional optimization of 
the perturbation size in an outer iteration, and we discuss algorithmic aspects.
In Section \ref{sect:eig-deriv} we provide the necessary theoretical background on first order
perturbation theory of eigenvalues and eigenvectors.
In Section \ref{sec:compr} we describe the complete computational procedure to compute the radius of
robustness.
In Section~\ref{sec:numexp} we present the results of numerical experiments where the 
robustness measure discussed here is computed for some examples of graphs.
Finally, in Section \ref{sec:concl} we draw our conclusions. 

\section{Related work}
\label{sec:related}
Related problems have been considered in the literature by several authors.  Most relevant to
our work are the papers \cite{CDM24,GVD24} and \cite{Nicosia}; see also the references therein.

In \cite{GVD24} the authors consider the problem of perturbing a given stochastic matrix so that
it has a prescribed stationary probability distribution, with the perturbation required to be minimal
with respect to the matrix 1-norm and to preserve the sparsity pattern of the original matrix (besides
of course, being such that the perturbed matrix remains stochastic). The solution is obtained by
formulating the problem as a constrained linear optimization problem which is then solved by
a suitable algorithm.

In \cite{CDM24}, the authors consider two popular measures of centrality, Katz and PageRank
(see, respectively, \cite{New10} and \cite{LM96}). Their goal is to find the ``minimal" perturbation
(according to a certain matrix functional) of the weights that enforce a prescribed vector of
centrality scores, while at the same time preserving the sparsity pattern.  In particular, the problem
of coalescing the top $m$ entries could be treated by the approach in \cite{CDM24}, by enforcing
a ranking vector in which the top $m$ entries take the same, prescribed value.  On the other hand,
PageRank and especially Katz centrality are not exactly the same as eigenvector centrality and,
perhaps more important, the optimality measure used in \cite{CDM24} is not the Frobenius norm
but a weighted sum of the (squared) Frobenius norm and the matrix 1-norm. Also, the approach
in \cite{CDM24}, based on quadratic programming, is quite different from the one we propose.
Nevertheless, also in view of the close relationship existing between Katz and eigenvector centrality
(see \cite{BeKl15}), it is possible that the approach in \cite{CDM24} may be applicable, with some
modifications, to the solution of the problem considered in the present paper.

Finally, we mention the interesting paper \cite{Nicosia}, which studies the broader problem of 
identifying sets of nodes in a network that can be used to ``control" the eigenvector centrality scores of 
the remaining nodes by manipulating the weights of the corresponding outgoing links. In this paper
it is shown how to identify minimal sets of nodes, called ``controller nodes," which can be used to 
obtain the desired eigenvector centrality ranking. This is clearly a different problem than the one
considered in this paper, where we seek the perturbation of minimum Frobenius norm, but the
approach in \cite{Nicosia} may be combined with the one in our paper, since our method allows
for the possibility of constraining the perturbations to a subset of the edges in the network
(see the example in Section \ref{sec:subp}).

\section{Robustness and vulnerability of eigenvector centrality}
\label{sec:problem}


Let 
\begin{equation*}
\mathcal{G}=(\mathcal{V},\mathcal E, A) 
\end{equation*}
be a weighted graph (either directed or undirected) 
with vertex set $\mathcal{V}=\{v_1,\dots,v_n \}$, edge set $\mathcal{E}\subset \mathcal{V}\times\mathcal{V}$
and weighted adjacency matrix $A$. {\color{black}
To each edge $(i,j)\in\mathcal{E}$ (from node $v_j$ to node $v_i$:  note the order) is associated a positive weight 
$a_{ij} >0 $; of course, for an undirected graph we have $a_{ji}=a_{ij}$ for all $i,j$.
We set $a_{ij}=0$ for $(i,j)\notin\mathcal{E}$. 
The weighted adjacency matrix of the graph is the non-negative matrix 
\[
A=(a_{ij}) \in \R^{n\times n}.
\]
When 
$G$ is undirected, $A$ is symmetric.}
To avoid a dependence of the measure on the scaling of the weights, we assume that $A$ is normalized with respect to the Frobenius norm,
\begin{equation} \label{eq:normalize}
\| A \| = \| A \|_F = 1
\end{equation}
where in this article - unless specified differently - $\| \cdot \| = \| \cdot \|_F$.

Given a graph $\widetilde{\mathcal{G}}=(\mathcal{V},\mathcal E, \widetilde A)$, with the same vertices and
edges of $\mathcal G$, we define the distance between $\mathcal G$ and $\widetilde{\mathcal{G}}$ by
\begin{equation} \label{eq:dist}
{\rm dist} \bigl( \mathcal G, \widetilde{\mathcal{G}} \bigr) = 
\sqrt{\sum\limits_{(i,j) \in \mathcal E} \left( a_{ij} - \widetilde a_{ij} \right)^2}
= \| A - \widetilde A \|\,.
\end{equation}

We order the eigenvalues $\{ \lambda_i \}$ of $A$ by non-increasing modulus so that $\lambda_1$ is the Perron
eigenvalue of $A$, which is non-negative and equal to $\rho(A)$ (the spectral radius of $A$).

\begin{assumption} \label{ass:sc}
In this article we assume that the graph is strongly connected, which is equivalent to the irreducibility of $A$.
\end{assumption}

The following theorem is classic (see e.g. \cite{HorJ91}).
\begin{theorem}[Perron--Frobenius] \label{th:PF}
If $A\in\R^{n\times n}$ is nonnegative and irreducible then
\begin{itemize}
\item[\rm (i) ] $\rho(A)>0$ is a simple eigenvalue of $A$;
\item[\rm (ii) ] there is a positive normalized eigenvector $\vp$ such 
that $A \vp = \rho(A) \vp$ ($\vp$ is called the Perron eigenvector of $A$).
\end{itemize}
\end{theorem}

Eigenvector centrality (see, e.g., \cite{New10}) measures the importance of the nodes in $\G$ by the magnitude of the entries of the  Perron eigenvector $\vp$,
allowing to rank the nodes. Note that the ranking is independent of the normalization used for $\vp$.

We indicate by $i_1, i_2, \ldots$ the indices of the entries of the Perron eigenvector $\vp$ of $A$, ordered by non-decreasing magnitude. 
 In the event that all the entries of $\vp$ are distinct,
the index $i_{1}$ corresponding to the largest entry in $\vp$ identifies the most central node $i_1$ in $\mathcal G$, as measured by eigenvector 
centrality; $i_2$ identifies the second most central node, and so forth.  Generally speaking, however, there could be many nodes with the
same centrality score. An extreme case is when all the row sums of $A$ are equal, in which case eigenvector centrality is unable to 
provide a ranking of the nodes, every node being given the same centrality score.

%
%
%
%
%
As an illustrative example, for the following graph $\G$ on the left picture,
\begin{center}
 \begin{tikzpicture}[
            > = stealth, 
            shorten > = 1pt, 
            auto,
            node distance = 3.3cm, 
            semithick, 
						main/.style = {draw, circle}
        ]

        \tikzstyle{every state}=[
            draw = black,
            thick,
            fill = white,
            minimum size = 4mm
        ]

        \node[main] (v0) {${\bf \emph{v}_1}$};
        \node[main] (v1) [above right of=v0] {$v_2$};
        \node[main] (v2) [right of=v0] {$v_3$};
        \node[main] (v3) [below right of=v0] {$v_4$};%

        \path[->] (v1) edge node {14.9} (v0);
        \path[->] (v2) edge node {6.7} (v0);
        \path[->] (v3) edge node {3.9} (v0);
				\path[->] (v0) edge node {4.8} (v3);
        \path[->] (v1) edge node {5.2} (v2);
        \path[->] (v2) edge node {4.0} (v1);
        \path[->] (v2) edge node {9.6} (v3);
				\path[->] (v3) edge node {14.1} (v2);

         0   14.9000    6.7000    3.9000
         0         0    4.0000         0
         0    5.2000         0   14.1000
    4.8000         0    9.6000         0
    \end{tikzpicture} 
		\hskip 2cm
		\begin{tikzpicture}[
            > = stealth, 
            shorten > = 1pt, 
            auto,
            node distance = 3.3cm, 
            semithick, 
						main/.style = {draw, circle}
        ]

        \tikzstyle{every state}=[
            draw = black,
            thick,
            fill = white,
            minimum size = 4mm
        ]

        \node[main] (v0) {$v_1$};
        \node[main] (v1) [above right of=v0] {$v_2$};
        \node[main] (v2) [right of=v0] {${\bf \emph{v}_3}$};
        \node[main] (v3) [below right of=v0] {$v_4$};%

        \path[->] (v1) edge node {14.9} (v0);
        \path[->] (v2) edge node {6.8} (v0);
        \path[->] (v3) edge node {4.0} (v0);
				\path[->] (v0) edge node {4.7} (v3);
        \path[->] (v1) edge node {5.1} (v2);
        \path[->] (v2) edge node {4.1} (v1);
        \path[->] (v2) edge node {9.5} (v3);
				\path[->] (v3) edge node {13.9} (v2);

    \end{tikzpicture}
\end{center}
we have 
\[
A = 
\left( \begin{array}{ccccc}
    0  &  14.9  &   6.7  &   3.9 \\
    0  &     0  &   4.0  &     0 \\
    0  &   5.2  &     0  &  14.1 \\
  4.8  &     0  &   9.6  &     0    
\end{array} \right), \quad 
\vp = 
\left( \begin{array}{c}
    0.5665 \\
    0.1570 \\
    \mathbf{0.5844} \\
    0.5594 
\end{array} \right). \hskip 0.1cm \ 
\]
{\color{black}
Hence, when ranked by (right) eigenvector centrality,
$v_3$ is the most central node,  $v_1$ the second most central one, and so forth. 

In this case, the ranking provided by eigenvector centrality appears unambiguous. 
However, a relative perturbation of the weights of the graph of 
a few percents causes the two
entries $\vp_1$ and $\vp_3$  to flip in magnitude, indicating a lack of robustness of the 
ranking.

We show this exhibiting by the nearby graph $\widetilde \G$ (on the right picture)  with weighted adjacency 
matrix $\widetilde A$ and
corresponding Perron eigenvector $\widetilde \vp$:
\[
\widetilde A = 
\left( \begin{array}{rrrr}
         0  &    14.9  &  6.8  &     4.0 \\
         0  &       0  &  4.1  &       0 \\
         0  &     5.1  &    0  &    13.9 \\
       4.7  &       0  &  9.5  &       0
\end{array} \right), \qquad
\widetilde \vp = 
\left( \begin{array}{c}
		\mathbf{0.5774} \\
    0.1602 \\
    0.5772 \\
    0.5548
\end{array} \right).
\]
Considering that in many applications the weights are known to an accuracy of only one or two digits, this example shows that
indeed the ranking determined by eigenvector centrality may not always be reliable.
}

Generally speaking, the ranking of the vertices obtained by this approach becomes unreliable when one of the following occurrences is 
determined by a small perturbation of the weights:

\begin{itemize}

\item[(I) ] The coalescence of the entries $\vp_{i_{1}}, \vp_{i_{2}}, \ldots, \vp_{i_m}$ 
for some index $m > 1$. 

\item[(II) ]  The coalescence of the eigenvalues $\lambda_1$ and $\lambda_2$.

\end{itemize}
Note that (II) can be ruled out if we make Assumption \ref{ass:sc}, i.e. in the irreducible case.

Indeed, if a small perturbation is able to coalesce the highest $m$ 
entries of the eigenvector $\vp$, it means that the ranking for the highest positions is not robust, especially in the presence of 
uncertainties, rounding errors,  or any other form of noise in the weights; 
in general this is a symptom of the fact that the highest $m$ vertices are difficult to rank using eigenvector centrality 
and it may therefore be better
to resort to different centrality measures.
%
A similar issue arises when a small perturbation makes the adjacency matrix $A$ reducible (i.e., the graph is no longer strongly
connected), in which case the Perron eigenvector may no longer be defined, rendering eigenvector centrality no longer usable. 
Note that issue (II) above falls under this scenario.

In this paper we direct our attention to issue (I), and  we are interested in computing the smallest perturbation of $A$ such that the
colaescence of the highest entries in the Perron vector occurs.
To this end, in the next subsection we propose a sharp measure which is aimed at characterizing the reliability of eigenvector centrality. 
As a by-product, the knowledge of this perturbation allows the possibility of modifying the ranking by small changes to the weights and 
may be used for this purpose.  Issue (II) is left as a topic for future work. 

\subsection{Robustness of the ranking and nearest ambiguous graph}
\label{sec:robustness measure}

We consider the general case of a directed graph and 
look for the closest strongly connected graph $\widetilde{\mathcal{G}} = (\mathcal{V},\mathcal E, \widetilde A)$ to 
$\mathcal G = (\mathcal{V},\mathcal E, A)$, such that the $m$ highest entries of the Perron eigenvector
coalesce. In order to preserve the strong connectivity we have to impose that the perturbed graph preserves the property 
$\widetilde a_{ij} > 0$ for $(i,j)\in\mathcal{E}$. 
For computational purposes we replace it with the property
\begin{equation} \label{eq:lbw}
\widetilde a_{ij} \ge \delta \qquad \mbox{for} \quad (i,j)\in\mathcal{E},
\end{equation}
with $1 \gg \delta > 0$ a small positive tolerance. 

We consider the following measure of robustness, which we call {\em robustness radius},
\begin{align*}
\varrho_m (\G) =\: \inf_{\widetilde{\G}} {\rm dist} \left( \mathcal G, \widetilde{\mathcal{G}} \right)
\stackrel{\triangle}=
\inf_{\widetilde A} \| A - \widetilde A\|,  
\end{align*} 
under the condition that the $m$ largest entries of the Perron eigenvector of $\widetilde A$ coalesce
and under the constraints $\widetilde a_{ij} \ge \delta$ for $(i,j)\in\mathcal{E}$
and  $\widetilde a_{ij}=0$ for $(i,j)\notin \mathcal{E}$.

Note that we may choose to only perturb a subset of the edges, $\mathcal{E}' \subset \mathcal{E}$, which would just
change the constraint to $\widetilde a_{ij} \ge \delta$ for $(i,j)\in\mathcal{E}'$. For example, we may choose to
perturb only the edges which either are subject to uncertainties or to changes in time.  Another possibility
would be to perturb only the nodes adjacent to a set of control nodes as considered, e.g., in \cite{Nicosia}.

Indeed it is possible (although highly unlikely) that no graph $\widetilde {\mathcal G}$ exists which fulfills the constraints, in which case
we say that the distance is $\infty$. Moreover, when there exist graphs fulfilling the constraints, we have that
the infimum is indeed a minimum, by compactness.  
\smallskip

\subsection{The considered problem.} We are interested in the following problem of eigenvector-constrained minimization.
\label{sec:consprob}
We let
\[
\Sigma = \left\{ B \in \R^{n\times n} \ \text{such that} \ 
b_{ij} \ge \delta \ \forall \ (i,j)\in\mathcal{E} \ \land \ 
b_{ij}=0 \ \forall \ (i,j)\notin \mathcal{E} \right\}
\]
and (denoting with $\up$ the Perron eigenvector of $B \in \Sigma$)  
\begin{equation} \label{eq:Omega}
\Omega = \left\{ B \in \Sigma : \up_{i_1} = \cdots = \up_{i_m} \right\}
\end{equation}
where, for a given $m$, $i_1, i_2, \ldots, i_m$ are the indices of the $m$ largest entries of the Perron eigenvector $\up$ of $B$
(i.e., $B \up = \rho(B) \up$).
The problem then reads
\begin{equation} \label{eq:nearness}
\min \| B - A \|_F \quad  {\rm s.\,t.} \quad B\in \Omega \,.
\end{equation}
 As we will see in the section on numerical experiments,  the constraint defining the set $\Omega$ is challenging to deal with by means of
classical optimization approaches.

In the next section we describe a reformulation of the problem which leads in
natural way to a two-level approach, an inner eigenvector optimization problem
and an outer root-finding one.

\subsection{Outline of the computational approach}
\label{sec:two-level method}
\label{sec:twolev}

In this section we describe an approach to compute the robustness measure $\varrho_m (\G)$ and the closest perturbed graph
characterized by an ambiguous ranking of the most important $m$ vertices. 
It is convenient to write 
\[
\widetilde A = A + \eps E \quad \mbox{ with} \ \| E \| = 1
\]
so that $\eps$ is the norm of the perturbation to $A$.

Based on the fact that we are addressing the coalescence of $m$ distinct entries, we can make use of their
average to obtain a concentration measure. For a given $\eps > 0$ and a given $E \in \R^{n \times n}$ with the
same pattern of $A$ (and such that $A + \eps E \ge 0$), we indicate by 
\begin{equation} \nonumber
\vp = \vp(A + \eps E)
\end{equation} 
the  (normalized)  Perron  eigenvector of $\widetilde A$ and define the mean of its largest $m$ entries as
\begin{equation} \label{eq:meanvp}
\Bigl\langle \vp(A+\eps E) \Bigr\rangle_m = \frac1m \sum_{j=1}^{m} \vp_{i_{j}} \bigl( A+\eps E \bigr),
\end{equation}
where $i_1, i_2, \ldots, i_m$ the indices of the $m$ largest entries of the Perron eigenvector $\vp$ of $A$, ordered by non-decreasing magnitude.

Then we define the following functional of $E$ :
\begin{eqnarray}\label{F-eps2}
F_\eps(E) & = & \frac12 \sum_{k=1}^{m} \Bigl( \vp_{i_{k}}  \bigl( A+\eps E \bigr) - \Bigl\langle \vp(A + \eps E) \Bigr\rangle_m \Bigr)^2
\end{eqnarray}
under the constraints of preservation of the sparsity pattern of~$A$, of  unit Frobenius norm of $E$ and of non-negativity $A+\eps E\ge 0$.
This nonnegative functional measures the  dispersion   of the largest $m$ entries of the Perron eigenvector
$\vp(A + \eps E)$ with respect to their mean. Clearly the functional vanishes if and only if the $m$ highest
entries of $\vp(A + \eps E)$ coalesce.

\smallskip

Our approach is summarized by the following two-level method:
\begin{itemize}
\item {\bf Inner iteration:\/} Given $\eps>0$, we compute a  smooth matrix flow 
$E(t)=(e_{ij}(t))\in\R^{n\times n}$, $t \ge 0$,  with the same sparsity pattern 
as $A$ (i.e., $e_{ij}(t)= 0$ if $a_{ij}=0$), with each $E(t)$ of unit Frobenius norm,  with 
$(A+\eps E(t))_{ij} \ge \delta$ (for $(i,j) \in \mathcal{E}$)  and we minimize
the functional $F_\eps(E(t))$ over this matrix flow.
The obtained minimizer is denoted by $E(\eps)$.

\item {\bf Outer iteration:\/} We compute the smallest value of $\eps$ such that the $m$-th largest entries of the Perron eigenvector 
$\vp \left( A+\eps E(\eps) \right)$ coalesce.
\end{itemize}

When there exists a solution to the problem $F_\eps \left( E(\eps) \right) = 0$,
in the outer iteration we compute the optimal $\eps$, denoted $\epstildem$
and the corresponding extremal perturbation $E(\epstildem)$, 
\begin{equation} \label{eq:epsm}
\epstildem \longrightarrow \arg\min\limits_{\eps > 0} F_\eps \left( E(\eps) \right) = 0
\end{equation}
by a combined Newton-bisection method. 

Although in most interesting cases the equation $F_\eps \left( E(\eps) \right) = 0$ has solutions, 
in general this is not assured.
If a solution to $F_\eps \left( E(\eps) \right) = 0$ exists then \eqref{eq:epsm} yields 
\[
\varrho_m (\G) = \epstildem.
\]

\subsubsection*{Existence of local minima}

An important issue concerns the existence of minima of \eqref{F-eps2}.
Let $\eps > 0$ be arbitrary but fixed.
By the assumption of strong connectivity of the perturbed graph, we have that the eigenvalue $\lambda = \rho(A + \eps E )$ is
simple for all admissible perturbations $E$. This implies the continuity with respect to $E$ of the Perron eigenvector $v$ 
associated to $\lambda$, which implies the continuity of the functional $F_\eps(E)$.
Combining this and the fact that the set of considered perturbations is compact, a minimizer of \eqref{eq:epsm} always exists. 
%
Since we make use of a gradient system approach to the problem, in general we cannot be sure to compute global optima, but only
to globally converge to local optima, which means that the final computed matrix provides an upper bound to the solution of  
\eqref{eq:epsm}.

\subsubsection*{Two-level optimization}
Due to the nested, or two-level form of problem \eqref{eq:Omega}--\eqref{eq:nearness}, it is difficult to devise a single-level 
optimization approach to its solution.

The particular two-level methodology we propose has two main motivations. The first is that in this way
we are able to approach the solution of the problem in a (norm) controlled way, which might be considered in analogy
to a trust region method, allowing us -- for example -- to reach the solution by a sequence of matrices of smaller norm.
The second is that in this way we obtain, as an important by-product, that the stationary points of the inner iteration
have an explicit characterization in terms of the constrained gradient of the
functional, therefore we are able to compute them explicitly by means of first-order eigenvector perturbation theory,
as shown in the  following Theorem \ref{thm:stat}. 
In particular, this reveals an underlying low-rank property that would otherwise remain hidden, and which can be exploited at 
the computational level.

For the solution of the inner iteration we propose a numerical integrator for the constrained gradient system associated
to the functional. 
This is based on an explicit formula for the gradient, which we derive by using first-order eigenvector perturbation theory.
Using alternative descent methods based on the exact knowledge of the gradient is also possible, and leads to similar
results.

For the solution of the outer iteration (which aims to solve a root-finding problem) we consider a Newton-bisection method, which exploits an exact
formula for the derivative of the function of which we aim to compute the smallest root.

\section{Derivatives of eigenvalues and eigenvectors}
\label{sect:eig-deriv}

In this section we recall a few basic results concerning first order perturbation
theory of real eigenvalues and eigenvectors of a real matrix, which is the relevant
case for this paper. 

\subsection{First order perturbation theory for simple eigenvalues}

The following standard perturbation result for eigenvalues is well-known; see  e.g. 
\cite[Theorem~1]{GreLO20} and \cite[Section II.1.1]{Kat95}.

\begin{theorem}[Derivative of simple eigenvalues]\label{lem:eigderiv} 
Consider a continuously differentiable path of square real matrices $M(t)$ for $t$ 
in an open interval $I$. Let $\lambda(t)$, $t\in I$, be a continuous path of real simple 
eigenvalues of $M(t)$. Let $x(t)$ and $y(t)$ be left and right eigenvectors, respectively, 
of $M(t)$ corresponding to the eigenvalue $\lambda(t)$. Then, $x(t)^\top y(t) \neq 0$ for $t\in I$ and 
$\lambda$ is continuously differentiable on $I$ with the derivative (we denote $\dot{\phantom a}=d/dt$)
\begin{equation}
\dot{\lambda} = \frac{x^\top \dot{M} y}{x^\top y}\,.
\end{equation}
Moreover, ``continuously differentiable'' can be replaced with ``analytic'' in the assumption and the conclusion.
\end{theorem}

Since we have  $x(t)^\top y(t) \neq 0$, we can apply the  normalization
\begin{equation}\label{eq:scalxy}
\| x(t) \|=1, \ \ \|y(t) \| =1, \quad x(t)^\top y(t) \text{ is real and positive.}
\end{equation}


\subsection{First order perturbation theory for eigenvectors}

For the derivative of eigenvectors we need the notion of  group inverse (or reduced resolvent); 
see \cite{MeyS88} as well as Kato (\cite[Section I.5.3]{Kat95}).

\begin{defi}[Group inverse]
\label{def:groupinv}
Let $M\in {\mathbb R}^{n\times n}$ be a singular matrix with a simple zero eigenvalue.
The \emph{group inverse} (or reduced resolvent) of $M$ is the unique matrix
$\Ginv$ with
\begin{equation}\label{group-inv-cond}
M \Ginv = \Ginv M, \qquad \Ginv M  \Ginv= \Ginv, \quad \mbox{and} \quad M \Ginv M = M.
\end{equation}
\end{defi}
It is known from Meyer and Stewart (\cite{MeyS88}) that if $M$ is normal, then
its group inverse $\Ginv$ is equal to the better known Moore--Penrose pseudoinverse 
$M^\dagger$. In general, the two pseudoinverses are not the same. 
They are, however, related by the following result, which is a special case of a more 
general result in \cite[Appendix A]{GugOS15} which is also verified directly.

\begin{theorem} [Group inverse and Moore--Penrose pseudoinverse]\label{thm:Ginv}
Suppose that the real matrix $M$ has the simple eigenvalue $0$ with corresponding left and 
right real eigenvectors  $x$ and $y$ of unit norm and such that $x^\top y > 0$.
Let $\Ginv$ be the group inverse of $M$, and with $\kappa=1/(x^\top y)$ define the projection 
$\Pi=\Id-\kappa x y^\top$.
The group inverse $\Ginv$ of $M$ is related to the Moore--Penrose pseudoinverse 
$M^\dagger$ by
\begin{equation}
\Ginv  =  \Pi M^\dagger \Pi .
\label{groupinvformula}
\end{equation}
\end{theorem}
This result is used in the codes implemented by the authors for the computations involving the group inverse.
Similar results are known in the Markov chain literature (see, e.g., \cite{Ben04}).
%
%
The group inverse appears in the following result on the derivative of eigenvectors, which is a variant 
of Theorem 2 of Meyer and Stewart (\cite{MeyS88}). 
\begin{theorem}[Derivative of eigenvectors] \label{thm:eigvecderiv}
Consider a smooth path of square real matrices $M(t)$ for $t$ in an open interval $I$. Let 
$\lambda(t)$, $t\in I$, be a path of simple real eigenvalues of $M(t)$. Then, there exists a continuously 
differentiable path of associated left and right eigenvectors $x(t)$ and $y(t)$, $t\in I$, which are 
of unit norm with $x^\top(t)y(t)>0$ and satisfy the differential equations
\begin{equation}\label{deigvec}
\begin{array}{rcl}
\dot x^\top\!\!\! & = & {}- x^\top \dot M \Ginv + \left(x^\top \dot M \Ginv x\right) x^\top,
\\[3mm]
\dot y  & = &  {}- \Ginv \dot M y  +  \left(y^\top \Ginv \dot M y\right) y,
\end{array}
\end{equation}
where $\Ginv(t)$ is the group inverse of $M(t)-\lambda(t) I$.

If $M$ is symmetric the previous simplify to $x=y$ and
\begin{equation}\label{deigvecsym}
\dot x^\top\!\!\! = {}- x^\top \dot M M^\dagger 
\end{equation}
where $M^\dagger(t)$ is the pseudoinverse of $M(t)-\lambda(t) I$.
\end{theorem}

Formulas \eqref{deigvec} are interpreted as follows: the first term in the right-hand side represents 
the free derivative of the eigenvector, while the second term constrains the eigenvector to preserve the
unit Euclidean norm with respect to $t$.

\section{A methodology for computing the closest ambiguous graph and the robustness radius}
\label{sec:compr}

In this section we describe the computational procedure to compute the robustness radius $\varrho_m(\G)$.
We indicate by 
\begin{equation*}
\langle X,Y \rangle = \mathrm{trace}(X^\top Y) 
\end{equation*}
the Frobenius inner product, i.e. the inner product inducing the Frobenius norm $\|\cdot\|=\|\cdot\|$  on $\R^{n\times n}$.

\subsection{Inner iteration}
We denote by $i_1, i_2, \ldots, i_m$ the indices of the $m$ largest entries of the Perron eigenvector $\vp \left( A + \eps E \right)$, 
ordered by non-decreasing magnitude. . Moreover, we denote the $i$-th unit versor by 
\[
\e_i = \Bigl( \underbrace{0 \ldots 0}_{i-1} \ 1 \ 0 \ldots 0 \Bigr)^\top  
\]
For a set of edges $ \mathcal{E}'$, we define $P_{\mathcal{E}'}$ as the orthogonal projection from $\R^{n\times n}$ onto the 
vector subspace of all real matrices having the sparsity pattern determined by~$\mathcal{E}'$: for $A=(a_{ij})$,
\[ 
P_{ \mathcal{E}'}(A)\big|_{ij} := \begin{cases}
	a_{ij}\,, &\text{if } (i,j)\in \mathcal{E}' \,,\\
	0\,, &\text{otherwise.}\end{cases}
\]


\subsection*{Further constraints}

An extension to the case where only some of the weights are subject to perturbations is straightforward.
For that, it is sufficient to define the subset $\mathcal{E}' \subset \mathcal{E}$ of those arcs whose
weights can be perturbed, and to replace everywhere the projection $P_{\mathcal{E}}$ by the analogous
projection $P_{\mathcal{E}'}$. 

Similarly, if a few weights are related by a linear relationship, this can be easily included in the 
considered projection with respect to the Frobenius inner product.   It is sufficient to determine
the projection onto the considered linear subspace and combine it with the one that takes into account
the sparsity structure, which is also expressed as a linear subspace constraint. For details, see subsection
\ref{sec:lincon}. 
%

From now on, unless otherwise specified, we shall consider $\mathcal{E}' = \mathcal{E}$, that is, we perturb
all edges of the graph.

\begin{defi}\label{def:adm}
For a given weight matrix $A$ and $\eps>0$, a matrix $E=(e_{ij})\in\R^{n\times n}$ is {\rm $\eps$-feasible} 
if 
\begin{itemize}
\item[(a)] $E$ is of unit Frobenius norm.
\item[(b)] $E=P_\E(E)$.   
\item[(c)] $(A+\eps E)_{ij} \ge \delta$ for $(i,j) \in \mathcal{E}$.  
\item[(d)] $E$ is symmetric (only for undirected graphs).
\end{itemize}
\end{defi}

\subsection{Structured gradient of the functional (\ref{F-eps2})}

Here we compute the gradient of the functional (\ref{F-eps2}) when restricted to matrices having the same sparsity pattern of $A$. 
The key idea is based on the fact that we are able to express the derivative of  $F_\eps(E)$ as $g(\dot E)$ for a special function $g$.
Using the simple fact that for vectors $b,c \in \R^n$ and a matrix $Z \in \R^{n \times n}$ one has
\[
b^\top Z c = \langle b c^\top, Z \rangle,
\] 
we are able to write (where $\dot E = \frac{d}{dt} E(t)$)
\[
g(\dot E) = \langle Z_\eps(E),\dot{E} \rangle 
\] 
with an appropriate matrix-valued function~$Z_\eps$, which -- in a natural way -- is interpreted as the constrained gradient of $F_\eps$ 
(with respect to the Frobenius inner product). 
Projecting it onto  the set of real matrices with sparsity pattern induced by $\E$  (first constraint) 
and unit norm (second constraint), we obtain the sought fully constrained gradient.

To address the undirected case, we define for $A \in \R^{n \times n}$,
\[ 
\sym(A)\big|_{ij} := 
	\frac{a_{ij} + a_{ji}}{2}
\]

Recalling \eqref{F-eps2}, we have that 
$F_\eps(E) = \frac12 \sum\limits_{k=1}^{m} \bigl( \vp_{i_{k}}  \left( A+\eps E \right) - \left\langle \vp(A + \eps E) \right\rangle_m \bigr)^2$.
After defining the auxiliary vector
\[
\bar\e = \frac1m \sum_{j=1}^{m} \e_{i_{j}}
\]
we are able to state the following result.

\begin{lemma}[Structured gradient of the functional \ -- \ directed graphs]  \label{lem:a-dot} 

Consider a smooth path $E(t)$ of $\eps$-feasible matrices, denote  by $\lambda(t) = \rho\left( A + \eps E(t) \right)$ the 
Perron eigenvalue and $\vp(t) = \vp \left( A + \eps E(t) \right)$ the associated Perron eigenvector.
Let $i_1, i_2, \ldots$ the indices of the entries of $\vp(t)$, ordered by non-decreasing magnitude.
We have (omitting the omnipresent dependence on $t$)
\begin{equation}\label{G-eps}
	\frac{d}{dt} F_\eps\left( E(t) \right) = \eta \Big\langle  G_\eps(E) , \dot{E} \Big\rangle, 
\end{equation}
where	$\eta > 0$ is a suitable scaling factor, 
with the structured gradient $G_\eps(E)$, having the sparsity pattern determined by the set of edges $\mathcal{E}$) 
\begin{equation} \label{eq:grad}
	G_\eps(E)=P_{\mathcal{E}} \left( R \right) 
\end{equation}
where $R$ is the free gradient:
\begin{eqnarray*}
R & = \sum\limits_{k=1}^{m} \left( \vp_{i_k} - \langle \vp\rangle_m \right)
	\Bigl( 
	\left( \vp_{i_k} - \langle \vp\rangle_m \right) (\Ainv)^\top \vp \vp^\top      
- (\Ainv)^\top \left( \e_{i_k} - \bar\e \right) \vp^\top \Bigr).
\end{eqnarray*}

\end{lemma}

\begin{proof}
By \eqref{F-eps2} we get
\begin{equation*}
\frac{d}{dt} F_\eps\left( E(t) \right) = \sum_{k=1}^{m}  \bigl(  \vp_{i_k} -  \langle \vp\rangle_m \bigr)
\bigl( \dot \vp_{i_k} -  \langle \dot \vp\rangle_m \bigr).
\end{equation*}

 Theorem \ref{thm:eigvecderiv} 
applied to $A + \eps E (t)$ yields 
\begin{eqnarray}
\dot \vp_{i_k} -  \langle \dot \vp\rangle_m & = & \Bigl( \e_{i_k} - \frac1m \sum_{j=1}^{m} \e_{i_{j}} \Bigr)^\top \dot \vp
= \eps \left( \e_{i_k} - \bar\e \right)^\top \left( - \Ainv \dot E \vp +  \left(\vp^\top \Ainv \dot E \vp\right) \vp\right)
\nonumber
\\
& = & \eps \Bigl( 
\big\langle - (\Ainv)^\top \left( \e_{i_k} - \bar\e \right) \vp^\top , \dot{E} \big\rangle
+ \left( \vp_{i_k} - \langle \vp\rangle_m \right) \big\langle (\Ainv)^\top \vp \vp^\top,  \dot E  \big\rangle 
\Bigr) 
\label{milestone1}
\end{eqnarray}
where $\Ainv(t)$ is the group inverse of $A + \eps E (t)-\lambda(t) I$ and $\vp(t)$ is the Perron eigenvector of
$A + \eps E(t)$, normalized with $\| \vp(t) \| = 1$.
Hence we can write
\begin{equation*}
\frac{d}{dt} F_\eps\left( E(t) \right) = \eps \big\langle  \Gf , \dot{E} \big\rangle
\end{equation*} 
with the free gradient 
\[
\Gf = \sum_{k=1}^{m} \bigl( \vp_{i_k} - \langle \vp\rangle_m \bigr)
	\Biggl( 
	\bigl( \vp_{i_k} - \langle \vp\rangle_m \bigr) (\Ainv)^\top \vp \vp^\top      
- (\Ainv)^\top \bigl( \e_{i_k} - \bar\e  \bigr) \vp^\top \Biggr).
\]

Using the property that $P_{\mathcal{E}} (\dot E) = \dot E$, we simply write
\begin{equation*}
\big\langle  \Gf , \dot{E} \big\rangle = \big\langle  P_{\mathcal{E}} (\Gf) + \left( \Gf - P_{\mathcal{E}} (\Gf) \right) , \dot{E} \big\rangle
\end{equation*}
and observe that 
\[
\big\langle  \Gf - P_{\mathcal{E}} (\Gf), \dot{E} \big\rangle = \big\langle  \Gf - P_{\mathcal{E}} (\Gf), P_{\mathcal{E}} (\dot E) \big\rangle = 0,
\]
This gives the structure-constrained gradient $G_\eps(E) = P_{\mathcal{E}} (\Gf)$ (see \eqref{eq:grad}).
\end{proof}

\begin{lemma}[Fully constrained gradient of the functional -- undirected graphs]  \label{lem:asym-dot} 
Under the same assumptions of Lemma \ref{lem:a-dot} to which is added the hypothesis that the smooth path $E(t)$ is also symmetric, 
we have that the fully constrained gradient is still given by \eqref{G-eps} with
\begin{equation}\label{eq:gradsym}
G_\eps(E) = {}-P_{\mathcal{E}} \left( \sum_{k=1}^{m} \left( \vp_{i_k} - \langle \vp\rangle_m \right)
	\sym \Bigl( A_\eps^\dagger \left( \e_{i_k} - \bar\e \right) \vp^\top \Bigr) \right) 
	\end{equation}
\end{lemma}

\begin{proof}
The proof is similar to that of Lemma \ref{lem:a-dot}. 
Here $\Ainv = A_\eps^\dagger$; then, 
by the well-known property $\vp^\top A_\eps^\dagger = 0$, 
we first get 
\[
\Gf = - (\Ainv)^\top \bigl( \e_{i_k} - \bar\e  \bigr) \vp^\top.
\]
Next, and including the symmetry constraint (ic) and thus the symmetry of $\dot E$, we get
\begin{equation*}
\big\langle  \Gf , \dot{E} \big\rangle = \big\langle  P_{\mathcal{E}} (\Gf) , \dot{E} \big\rangle =
\big\langle  \sym \left( P_{\mathcal{E}} (\Gf) \right) + \skw \left( P_{\mathcal{E}} (\Gf) \right) , \dot{E} \big\rangle
=
\big\langle  \sym \left( P_{\mathcal{E}} (\Gf) \right), \dot{E} \big\rangle
\end{equation*}
where for a matrix $B$, $\skw(B) = {(B - B^\top)}/{2}$. Orthogonality of symmetric and skew-symmetric matrices
proves the result, i.e. gradient \eqref{eq:grad} symplifies to \eqref{eq:gradsym}.
\end{proof}

\subsection{Constrained gradient of the functional}

In this section we take into account the norm constraint on $E$ and the inequality constraint on $A + \eps E$.
One possibility would be to use a penalization method; this, however, would have the disadvantage that the
potential violation of the inequality constraints might have an impact on the Perron eigenvector and on its 
continuity. For this reason we opt for a Lagrangian approach. 

\subsubsection*{Admissible directions}
We recall the considered constraints on $E$:
\begin{itemize}
\item[(a)] $E$ is of unit Frobenius norm.
\item[(b)] $E=P_\E(E)$.
\item[(c)] $(A+\eps E)_{ij} \ge \delta$ for $(i,j) \in \mathcal{E}$.
\item[(d)] $E$ is symmetric (only for undirected graphs).
\end{itemize}

Since $E(t)$ is of unit Frobenius norm by condition (a), we have
\[ 
0 = \frac12\, \frac{d}{dt} \lVert E(t) \rVert^2 =  \langle E(t), \dot E(t) \rangle .
\]
To fulfill the non-negativity condition (c) of the  matrix $A+\eps E(t)$,  we need that 
\[
\dot{e}_{ij}(t)  \ge 0 \quad  \mbox{for all} \ (i,j) \ \mbox{with} \
a_{ij} + \eps e_{ij}(t) = \delta \ge 0.
\]  

Conditions (b) and (d) are satisfied if the same condition 
holds for $\dot{E}$. These four conditions are in fact also sufficient for a matrix to be the time derivative of a path of $\eps$-feasible matrices.

Hence, for every $\eps$-feasible matrix $E$, a matrix $Z=(z_{ij})\in\R^{n\times n}$ is the derivative  of some path of $\eps$-feasible matrices passing through $E$ if and only if the following four conditions are satisfied:
\begin{itemize}
\item[(a')] $\langle E,Z\rangle=0$.
\item[(b')] $Z=P_\E(Z)$.
\item[(c')] $P_{\E_0}(Z)\ge 0$.
\item[(d')] $Z$ is symmetric (only for undirected graphs) .
\end{itemize}
Here  $\mathcal{E}_0=\E_0(\eps E)$ is the set 
\[ 
\mathcal{E}_0 := \{(i,j) \in \mathcal{E}:\, a_{ij} + \eps e_{ij} = \delta \}.   
\]
Condition (c') implies that $z_{ij} \ge 0 \text{ for all }(i,j) \in \mathcal{E}_0$.

\subsubsection*{Admissible direction of steepest descent}
In order to determine the admissible direction $Z = \dot E$ of steepest descent from $E$, 
we  therefore consider the following optimization problem for $G=G_\eps(E)$:
\begin{align}\label{opt_problem_1}
& \min_{Z} \langle G,Z \rangle \quad\text{subject to} \\
& \begin{cases}
		\langle E, Z \rangle = 0 \,, \\
		z_{ij} \ge 0 \text{ for all }(i,j) \in \mathcal{E}_0 \,,\\
		z_{ij} = z_{ji}  \text{ for all }(i,j) \quad \text{(for undirected graphs})\,,\\
		\langle Z,Z \rangle = 1 \,.
	\end{cases}
\end{align}
The additional constraint  $\lVert Z \rVert = 1$ just normalizes the descent direction. 
Problem \eqref{opt_problem_1} has a quadratic constraint. 

We now formulate a quadratic optimization problem with linear constraints, which is equivalent 
to \eqref{opt_problem_1}.

\subsection*{An equivalent quadratic optimization problem with linear constraints}
\begin{align}\label{opt_problem_2}
	& \min_{Z} \langle Z,Z \rangle \quad \text{subject to}  \\
	& \begin{cases}
			\langle E, Z \rangle = 0 \,, \\
			 z_{ij} \ge 0 \text{ for all }(i,j) \in \mathcal{E}_0 \,,\\
			 z_{ij} = z_{ji}  \text{ for all }(i,j) \quad \text{(for undirected graphs})\,,\\
			\langle G, Z \rangle = -1 \,.
		\end{cases}
\end{align}
Problems \eqref{opt_problem_1} and \eqref{opt_problem_2} are equivalent in the sense that Problem \eqref{opt_problem_2} 
yields the same solution, i.e. descent direction, provided that a strict descent direction exists, i.e., satisfying 
$\langle G,Z \rangle<0$ and the constraints (a')--(d'). 
The equivalence is a consequence of the fact that when  $\langle G, Z \rangle < 0$, there exists a scaling factor $\alpha>0$ s
uch that  $\langle G, \alpha Z \rangle = -1$.

Apart from the normalization,
problems \eqref{opt_problem_1} and \eqref{opt_problem_2} yield the same KKT (Karush--Kuhn--Tucker) conditions. 
Since the objective function $\langle Z,Z \rangle$ of problem \eqref{opt_problem_2} is convex and all constraints are 
linear, the KKT conditions are not only necessary but also sufficient conditions (\cite[Theorem 9.4.1]{F13}),  that is, 
a KKT point is already a solution of the optimization problem.

The solution of \eqref{opt_problem_2} satisfies the KKT conditions (see, e.g., \cite[Ch. 2.2]{geigerkanzow})
\begin{align}  \label{KKT}
	& \begin{cases} 
	\kappa Z = - G - \gamma E + \sum_{(i,j) \in \mathcal{E}_0} \mu_{ij} \e_i \e_j^T \,, \\
	z_{ij} \ge 0 \text{ for all }(i,j) \in \mathcal{E}_0 \,, \\
	\langle E, Z \rangle = 0 \,, \\
	\langle Z, Z \rangle = 1 \,, \\
	 Z = Z^\top \quad \text{(for undirected graphs})\,,\\
	\mu_{ij} z_{ij} = 0 \text{ for all }(i,j) \in \mathcal{E}_0 \,, \\
	\mu_{ij} \ge 0 \text{ for all }(i,j) \in \mathcal{E}_0 \,,
	\end{cases}
\end{align}
where $\e_\ell$ is the $\ell$-th unit versor, $\gamma$ is determined by 
the normalization $\langle E, Z \rangle = 0$, and $\kappa > 0$ is determined by 
the normalization $\lVert Z \rVert^2 = \langle Z, Z \rangle = 1$. 

\subsection{Constrained gradient flow}
The gradient flow of $F_\eps$ under the constraints (a)--(d) is given by
\begin{equation}
\label{eq:ode}
\dot E(t) = Z(t),
\end{equation}
where $Z(t)$ solves the KKT system \eqref{KKT} with $G=G_\eps(E(t))$ under the constraints (a')--(d'). 
Note that the set of edges $\E_0(t)=\E_0(\eps E(t))$. 

\begin{lemma} On an interval where $\E_0(t)$ does not change, the gradient system becomes, with $\Pact =P_{\E \setminus \E_0}$,
\begin{equation}\label{Pact-ode}
\dot E = - \Pact  G_\eps(E) + \gamma \Pact  E \quad\hbox{ with }\quad
\gamma = \frac{\langle G_\eps(E), \Pact  E \rangle }{\| \Pact  E\|^2}.
\end{equation}
\end{lemma}

\begin{proof} The positive Lagrange multipliers $\mu_{ij}>0$ just have the role to ensure that $\dot e_{ij}=0$.  
With $G=G_\eps(E)$, the gradient system therefore reads
\begin{equation}\label{E-ode}
\dot E = \Pact (-G + \gamma E),
\end{equation}
where $\gamma$ is determined from the constraint $\langle E,\dot E \rangle =0$. We then have
\begin{equation*}
0=\langle E,\dot E \rangle = \langle E,\Pact (-G + \gamma E) \rangle=
- \langle \Pact  E,G \rangle +  \gamma\langle \Pact  E,\Pact  E \rangle,
\end{equation*}
and the result follows.
\end{proof}

%
%
%

\subsection{Monotonicity and stationary points}

The following result shows the monotonicity of the functional  $F_\eps (E(t))$ with respect to $t$ and provides
a characterization of the stationary points. We observe here that monotonicity follows directly from the construction of the gradient system \eqref{Pact-ode} 
(for the time being, we assume that the non-negativity constraint does not activate so that penalization is not necessary, 
which is what we observe in all the tests we have done).

\begin{theorem} \label{thm:stat}  
Let $E(t)$ of unit Frobenius norm satisfy the differential equation \eqref{Pact-ode} with $G_\eps(E)$ as in \eqref{eq:grad}. 
Then, the functional  $F_\eps (E(t))$ decreases monotonically with respect to $t$:
\[
\dot F_\eps(E(t)) \le 0.
\]
Furthermore, the following statements are equivalent along solutions of \eqref{Pact-ode}:
\begin{enumerate}
\item ${\dot F}_{\eps} = 0$.
\item $\dot E =0$.
\item $\Pact E$ is a real multiple of $\Pact G_\eps(E)$.
\end{enumerate}
\end{theorem}

\begin{proof} 
We prove the result in three steps.
\begin{itemize}
\item[] (1) implies (3): Using Lemma~\ref{lem:a-dot} and \eqref{Pact-ode} we obtain, with $G=G_\eps(E)$,
\begin{equation} \label{F-monotone}
 \frac1\eps\, \dot F_\eps = \langle G, \dot E \rangle 
 = \langle G, - \Pact G + \gamma \Pact E \rangle
 = - \| \Pact G \|^2 + 
 \frac{\langle \Pact G, \Pact E \rangle^2}{\|\Pact E\|^2}.
\end{equation}
By means of the strong form of the Cauchy--Schwarz inequality we obtain monotonicity. Furthermore we have that 
\[
\frac{d}{dt} F_{\eps} \left( E(t) \right) = 0 \qquad \Longrightarrow \qquad \mbox{$\Pact E$ is a multiple of $\Pact G$}.
\] 

\item[] (3) implies (2): This is a direct consequence of the vanishing of the r.h.s. of \eqref{Pact-ode}.

\item[] (2) implies (1): This is also a direct consequence of Lemma \ref{lem:a-dot}.
\end{itemize}
\end{proof}


\subsection{Underlying low-rank property of extremizers} \label{sec:lowrank}

For any graph, we have a number of ODEs which is given by the nonzero entries  in the adjacency matrix  
associated to the graph. Hence, for a dense graph we may have up to $n^2$ ODEs.
If inequality constraints are inactive, the gradient system \eqref{Pact-ode} is such that $\Pact$ is the identity and 
the gradient system therefore reads
\begin{equation}\label{E-ode-no}
\dot E = -G + \gamma E, \qquad \gamma = \langle G, E \rangle.
\end{equation}
where $G=G_\eps(E)$.
In this case extremizers (stationary points) are characterized by an underlying low-rank property, being
$G$ the projection of a low-rank matrix onto the sparsity-pattern of the weighted adjacency matrix.
Clearly, if $\eps < \min_{ij} a_{ij}$, inequality constraints cannot be active so that we can ignore them
and solve the system of ODEs \eqref{E-ode-no}.
Moreover, in all our experiments we have never observed activation of inequality constraints associated to 
nonnegativity.

Following the lines of \cite{GLS23}, we are able to exploit the underlying low-rank structure of the
solutions of the optimization problem we consider.
By Theorem \ref{thm:stat} we have that the stationary points which we aim to compute are proportional to
the gradient $G$, which - as we have seen - is the projection of a rank-$1$ matrix. Proceeding similarly to
\cite{GLS23}, by writing $E = P_{\mathcal{E}}(H)$  for a rank-$1$ matrix $H$, 
we are able to obtain a new system of ODEs for $H$, whose dynamics is constrained to the manifold of rank-$1$ 
matrices.
We can prove that the stationary points of the two systems of ODEs are in a one-to-one correspondence
and also that the system of ODEs
for $H$, even if it loses the gradient system structure, converges - at least locally - to the stationary points.
We omit a full derivation for sake of conciseness and refer the reader to \cite{GLS23} for details.

\subsection{Including a linear constraint on the entries of the adjacency matrix}
\label{sec:lincon}

Suppose now that a few weights associated with the edge set $\mathcal{E}' \subseteq \mathcal{E}$
are related by a linear relationship, that is
\begin{equation} \label{eq:lincon}
\sum_{(i,j) \in \mathcal{E}} b_{ij} a_{ij} = c ,
\end{equation}
with $c, b_{ij} \in \R$ for all $(i,j) \in \mathcal{E}'$ and $b_{ij} =0$ for all $(i,j) \in \mathcal{E} \setminus \mathcal{E}'$. 
Then for the perturbed matrix $A + \eps E(t)$ we have
to impose for all $t$ 
\[
\sum_{(i,j) \in \mathcal{E}} b_{ij} e_{ij}(t) = 0 \quad \Longrightarrow \quad \sum_{(i,j) \in \mathcal{E}} b_{ij} \dot e_{ij}(t) = 0,
\]  
that is, we have to impose the scalar product $\langle b, \dot{e}(t) \rangle = 0$ for the associated vectors $b = \left( b_{ij} \right)$ and 
$e = \left( e_{ij} \right) = {\rm vec}(E)$, which represents a column vector formed by the nonzero entries of $E$ (which are defined by the 
sparsity pattern).

This implies that we have to modify the r.h.s.~in the matrix ODE \eqref{eq:ode} in the inner iteration, which we write in short as
$\dot{E} = Z$. Writing it in vectorized form as $\dot{e} = z$, with $z = {\rm vec}(Z)$, we modify it as
\begin{equation} \label{eq:odevec}
\dot{e} = z - \frac{1}{\| b \|^2}\langle z, b \rangle b ,
\end{equation}
that is, we orthogonally project the vector field $z$ onto the linear space ${\rm Ker}(b) = \langle b, \cdot \rangle = 0$.
Indeed, this implies
\[
\langle \dot{e}, b \rangle = \langle z, b \rangle - \frac{1}{\| b \|^2} \langle z, b \rangle \langle b, b \rangle  = 0 ,
\]
as desired; it follows that the property \eqref{eq:lincon} is conserved along the flow of the modified ODE.
Equation \eqref{eq:ode} is in vector form, but it can be rewritten directly in matrix form.

\section{Undirected graphs}

Undirected graphs lead to some significant simplifications.
In this case we have that $\Ainv = A_\eps^\dagger$, the Moore-Penrose pseudoinverse,
which is symmetric.

Recalling \eqref{eq:gradsym}, we have 
\begin{eqnarray} \nonumber 
	G_\eps(E) & = & {}-P_{\mathcal{E}} \left( \sum_{k=1}^{m} \left( \vp_{i_k} - \langle \vp\rangle_m \right)
	\sym \Bigl( A_\eps^\dagger \left( \e_{i_k} - \bar\e \right) \vp^\top \Bigr) \right) 
\\[2mm]	\nonumber
	& = &
	{}-P_{\mathcal{E}} \left( \sum_{k=1}^{m} c_k
	\sym \Bigl( a_k \vp^\top \Bigr)		
	\right),
\end{eqnarray}
where we have defined 
\[
a_k = A_\eps^\dagger b_k = (A + \eps E-\lambda I)^\dagger b_k, \qquad \mbox{with} \qquad 
b_k = \left( \e_{i_k} - \bar\e \right), \qquad
c_k = \bigl( \vp_{i_k} - \langle \vp\rangle_m \bigr),
\]
which may be computed by solving the nonsingular linear system
\begin{equation}\label{a-ls}
\begin{pmatrix}
A + \eps E-\lambda I & \vp\\
\vp^\top & 0
\end{pmatrix}
\begin{pmatrix}
a_k \\
\mu 
\end{pmatrix}
=
\begin{pmatrix}
b_k \\
0 
\end{pmatrix}.
\end{equation}
Note that the submatrix $A + \eps E-\lambda I$ above is negative semidefinite and singular. 
It is straightforward to check that  the linear system \eqref{a-ls} is nonsingular if and only if the orthogonal complement of 
$\vp$ and the kernel of $A + \epsilon E - \lambda I$ 
have trivial intersection.
 Since $\text{Ker} (A + \eps E-\lambda I) = \text{Span}\,\{\vp\}$ and obviously $\vp^\top \vp\ne 0$, the condition is satisfied and the system \eqref{a-ls}
has a unique solution. 


The computational cost of computing $G_\eps(E)$ consists of computing the Perron eigenvector of the matrix $A+\eps E(t)$ for $t$ a suitable point on the mesh of integration and in solving -- at every step -- the linear system \eqref{a-ls}. The Perron eigenvector can be computed efficiently, for example using Lanczos' method; if the Perron root is well
separated by the rest of the spectrum, even the power method will converge quickly. The solution of the linear system \eqref{a-ls} merits a more detailed discussion.

\subsubsection*{Solving the linear system \eqref{a-ls}}

In the considered undirected case, system \eqref{a-ls} is symmetric and indefinite, with one positive and $n$ negative eigenvalues.
For small and moderate values of $n$, the solution can be obtained by means of some variant of the LDL$^T$ factorization, see \cite{GVL4}.
For large and sparse problems, an efficient way to solve \eqref{a-ls} is by means of the Minimal Residual (MinRes) iteration \cite{PS75}, possibly combined 
with a suitable symmetric positive definite preconditioner \cite{BGL05}. Alternatively, system \eqref{a-ls} can be replaced by the equivalent one
\begin{equation}\label{a-ls-alt}
\begin{pmatrix}
A + \eps E-\lambda I - \vp \vp^\top & \vp \\
\vp^\top & 0
\end{pmatrix}
\begin{pmatrix}
a_k \\
\mu 
\end{pmatrix}
=
\begin{pmatrix}
b_k \\
0 
\end{pmatrix}.
\end{equation}
Note that the submatrix $A + \eps E-\lambda I - \vp \vp^\top$ is negative definite (nonsingular).  This matrix is completely full, since $\vp>0$,
and it should not be formed explicitly.  System \eqref{a-ls-alt} can be solved by means of the block LU factorization
\begin{equation*} 
{\small
\begin{pmatrix}
A + \eps E-\lambda I - \vp \vp^\top& \vp \\
\vp^\top & 0
\end{pmatrix}
=
\begin{pmatrix}
A + \eps E-\lambda I - \vp \vp^\top & 0 \\
\vp^\top  &  - s
\end{pmatrix}
\begin{pmatrix}
I  &  (A + \eps E-\lambda I - \vp \vp^\top)^{-1} \vp \\
0 & 1
\end{pmatrix}
}
\end{equation*}
where we have set $s =\vp^\top(A + \eps E-\lambda I - \vp \vp^\top)^{-1} \vp$. It follows that the solution of \eqref{a-ls-alt}, and
therefore of \eqref{a-ls}, is given by
\begin{eqnarray*}
a_k & = & (A + \eps E-\lambda I - \vp \vp^\top)^{-1}b_k - \mu (A + \eps E-\lambda I - \vp \vp^\top)^{-1} \vp\,,  
\\
\mu & = & 
\frac{\vp^\top(A+ \eps E-\lambda I - \vp \vp^\top)^{-1}b_k}{\vp^\top(A+ \eps E-\lambda I - \vp \vp^\top)^{-1} \vp}\,,
\end{eqnarray*}
showing that this approach requires the solution of two linear systems with the same coefficient matrix
$A + \eps E-\lambda I - \vp \vp^\top$ and right-hand sides $b_k$ and $\vp$, respectively.  These two linear systems can be
solved by means of the conjugate gradient (CG) method, possibly combined with a preconditioner. Clearly, the CG
 method is well suited since matrix-vector products involving $A + \eps E-\lambda I - \vp \vp^\top$ can be computed
without forming the coefficient matrix explicitly; for sparse $A$, the cost of each matrix-vector product is linear in $n$.
The choice of preconditioner is clearly problem dependent; note, however, that as the outer iteration converges, the coefficient
matrix of the linear systems varies slowly from one step to the next, making it possible to reuse the same preconditioner
for several steps.

\subsubsection*{Projected rank-$2$ properties of extremizers}

Assume the inequality constraints are inactive. We let
\[
L_\eps(E) = \left( \sum_{k=1}^{m} \left( \vp_{i_k} - \langle \vp\rangle_m \right)
	\sym \Bigl( A_\eps^\dagger \left( \e_{i_k} - \bar\e \right) \vp^\top \Bigr) \right)
\]
and $G_\eps(E) = {}-P_{\mathcal{E}} \left( L_\eps(E) \right)$.
Note that $A_\eps^\dagger \left( \e_{i_k} - \bar\e \right) \vp^\top$ is a rank-$1$ nilpotent matrix.
We let
\[
g_{i_j} = A_\eps^\dagger \e_{i_j}, \qquad \mbox{the $i_j$-th column of} \ A_\eps^\dagger
\]
and
\[
\tau_{i_k} = \vp_{i_k} - \langle \vp\rangle_m, \qquad \qquad a_{i_k} = g_{i_k} - \langle g \rangle_m, \quad \mbox{with} \quad
\langle g \rangle_m = \frac1m \sum_{j=1}^{m} g_{i_{j}}.   
\]
This notation allows us to rewrite
\[
L_\eps(E) = \sym \left( a \vp^\top \right), \qquad a = \sum_{k=1}^{m} \tau_{i_k} a_{i_k}, 
\]
which is a rank-$2$ symmetric matrix.

If $(i,j)\in\mathcal{E}$ for all $i=1,\ldots,n$ and $j=1,\ldots,n$ (which includes possible self-loops), we have that the pattern $\mathcal{E}$ 
is full, so that $P_{\mathcal E}$ is the identity projector. Thus $G_\eps(E) = - L_\eps(E)$, which means that the stationary points 
(if the penalty function is inactive) are rank-$2$.

In the general case, stationary points are projections onto $\mathcal{E}$ of rank-$2$ symmetric matrices.


\subsection{An illustrative example of an undirected graph}
\label{sec:ill}


\tikzstyle{vertex}=[circle,fill=black!10,minimum size=20pt,inner sep=0pt]
\tikzstyle{selected vertex} = [vertex, fill=red!24]
\tikzstyle{second selected vertex} = [vertex, fill=blue!24]
\tikzstyle{edge} = [draw,thick,-]
\tikzstyle{weight} = [font=\small]
\tikzstyle{selected edge} = [draw,line width=5pt,-,green!50]
\tikzstyle{ignored edge} =  [draw,line width=5pt,-,blue!20]

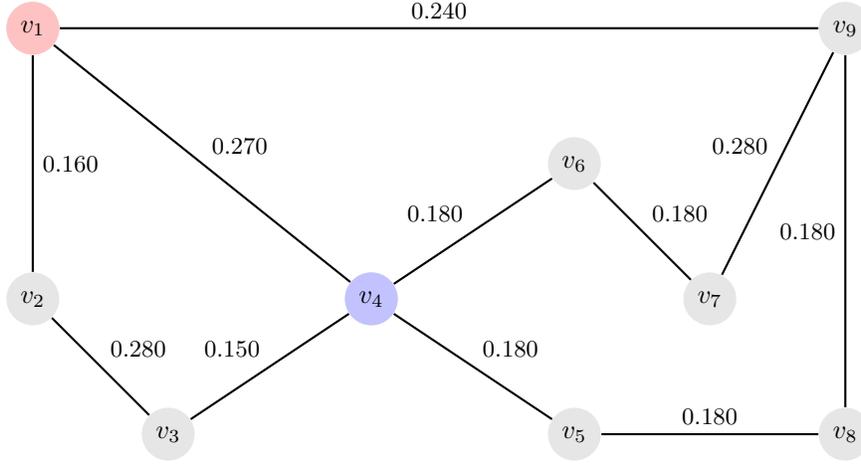
\begin{figure}[ht]
\begin{center}
\begin{tikzpicture}[scale=1.8, auto,swap]
    \foreach \pos/\name in {{(0,2)/v_1}, {(1,-1)/v_3}, {(4,1)/v_6}, {(6,2)/v_9},
                            {(0,0)/v_2}, {(4,-1)/v_5}, {(2.5,0)/v_4}, {(5,0)/v_7}, {(6,-1)/v_8}}
        \node[vertex] (\name) at \pos {$\name$};
    \foreach \source/ \dest /\weight in {v_2/v_1/0.160, v_4/v_1/0.270, v_9/v_1/0.240,
                                         v_3/v_2/0.280, 
																				 v_4/v_3/0.150,
																				 v_5/v_4/0.180,v_6/v_4/0.180,
                                         v_8/v_5/0.180,
																				 v_7/v_6/0.180,
                                         v_9/v_7/0.280,
																				 v_9/v_8/0.180}
        \path[edge] (\source) -- node[weight] {$\weight$} (\dest);
        \foreach \vertex / \fr in {v_1/1} 
				\node[selected vertex] at (\vertex) {$\vertex$};
				\foreach \vertex / \fr in {v_4/2} 
				\node[second selected vertex] at (\vertex) {$\vertex$};
\end{tikzpicture}
\end{center}
\caption{Example of symmetric weighted graph}
\label{fig:2}
\end{figure}

Consider the simple undirected graph with $9$ vertices of Figure \ref{fig:2}.
From a topological point of view one could argue that the vertex $v_4$ is characterized
by the highest centrality score; however, since the graph is weighted this is not
true in general.

In order to associate to the graph a scaled adjacency matrix $A$ with unit Frobenius norm,
we multiply the adjacency matrix whose weights are indicated in Figure \ref{fig:2} by the 
normalizing factor 
\[
\beta \approx 1/0.9974 \approx 1.0026. 
\]

We have $\lambda = \rho(A) \approx 0.5559\ldots$ and the associated Perron eigenvector (to $4$ digits)
\small
\begin{equation*}
\vp \approx \left( 
\begin{array}{ccccccccc}
    \mathbf{0.4844} &
    0.2712 &
    0.2602 &
    \mathbf{0.4553} &
    0.2154 &
    0.2433 &
    0.2941 &
    0.2082 &
    0.4259 
\end{array}
\right)^\top.
\end{equation*}
\normalsize

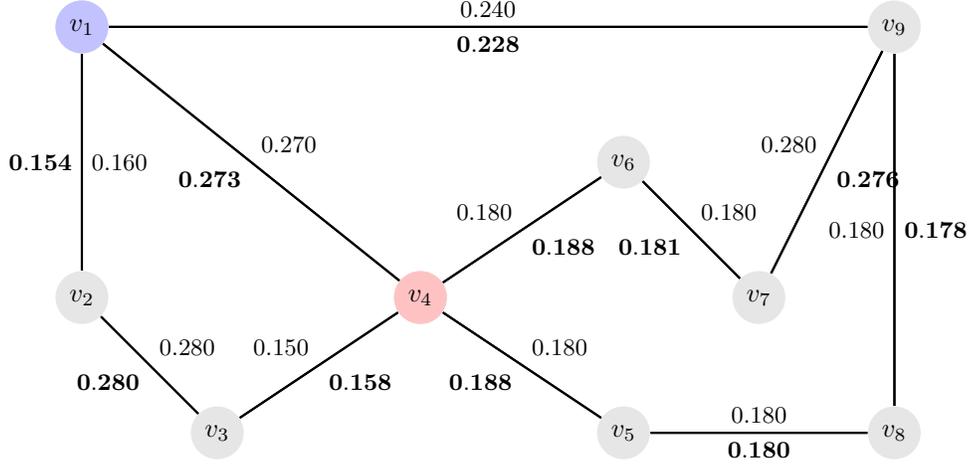
\begin{figure}[ht]
\begin{center}
\begin{tikzpicture}[scale=1.8, auto,swap]
    \foreach \pos/\name in {{(0,2)/v_1}, {(1,-1)/v_3}, {(4,1)/v_6}, {(6,2)/v_9},
                            {(0,0)/v_2}, {(4,-1)/v_5}, {(2.5,0)/v_4}, {(5,0)/v_7}, {(6,-1)/v_8}}
        \node[vertex] (\name) at \pos {$\name$};
    \foreach \source/ \dest /\weight in {v_2/v_1/0.160, v_4/v_1/0.270, v_9/v_1/0.240, 
                                         v_3/v_2/0.280, 
																				 v_4/v_3/0.150,
																				 v_5/v_4/0.180,v_6/v_4/0.180,
                                         v_8/v_5/0.180,
																				 v_7/v_6/0.180,
                                         v_9/v_7/0.280,
																				 v_9/v_8/0.180,
																				 v_1/v_2/\mathbf{0.154},v_1/v_4/\mathbf{0.273},v_1/v_9/\mathbf{0.228},
																				 v_2/v_3/\mathbf{0.280},
																				 v_3/v_4/\mathbf{0.158},
																				 v_4/v_5/\mathbf{0.188},v_4/v_6/\mathbf{0.188},
																				 v_5/v_8/\mathbf{0.180},
																				 v_6/v_7/\mathbf{0.181},
																				 v_7/v_9/\mathbf{0.276},
																				 v_8/v_9/\mathbf{0.178}}
        \path[edge] (\source) -- node[weight] {$\weight$} (\dest);
        \foreach \vertex / \fr in {v_4/2} 
				\node[selected vertex] at (\vertex) {$\vertex$};
				\foreach \vertex / \fr in {v_1/1} 
				\node[second selected vertex] at (\vertex) {$\vertex$};
\end{tikzpicture}
\end{center}
\caption{Perturbed graph of Figure \ref{fig:2} (weights of modified graph appear in bold on each edge).}
\label{fig:3}
\end{figure}

{\color{black} The ranking of nodes $\vp_1$ and $\vp_4$ can be interchanged by small changes in the weights. This
is illustrated in Figure \ref{fig:3}.}
Choosing a suitable small perturbation (of norm $0.0280$) of the normalized adjacency matrix $A$, 
we obtain the modified Perron eigenvector $\tilde \vp$ 
\small
\begin{equation*}
\tilde \vp = \left( 
\begin{array}{ccccccccc}
    \mathbf{0.4754} &
    0.2699 &
    0.2722 &
    \mathbf{0.4758} &
    0.2277 &
    0.2543 &
    0.2839 &
    0.2035 &
    0.4027
\end{array}
\right)^\top.
\end{equation*}
\normalsize
where now $\vp_4$ precedes $\vp_1$.

Note that, by a simple continuity argument of the Perron eigenvector, there exists a perturbation of the weights for which $\tilde \vp_1$ and $\tilde \vp_4$ have exactly the
same value.  

We are interested in the robustness of the ranking for the first $m=2$ positions, corresponding to 
$i_1 = 1$ and $i_2 = 4$.
The gap between the $1$-st and $4$-th entries is
$
\vp_1 - \vp_4 = 0.0291\ldots
$

According to \eqref{eq:gradsym}, we compute the structure-constrained gradient $G_0$, normalized to have Frobenius norm equal to $1$; it is
given by 
\small
\begin{equation*}
G_0 = \left( \begin{array}{rrrrrrrrr}
    0      &\!   0.336  &\!    0      &\!   -0.113 &\!    0      &\!    0      &\!    0      &\!    0      &\!    0.645  \\
    0.336  &\!    0     &\!   -0.029  &\!    0     &\!    0      &\!    0      &\!    0      &\!    0      &\!    0      \\
    0      &\!  -0.029  &\!    0      &\!   -0.412 &\!    0      &\!    0      &\!    0      &\!    0      &\!    0      \\
   -0.113  &\!   0      &\!   -0.412  &\!    0     &\!   -0.418  &\!   -0.438  &\!    0      &\!    0      &\!    0      \\
    0      &\!   0      &\!    0      &\!   -0.418 &\!    0      &\!    0      &\!    0      &\!   -0.073  &\!    0      \\
    0      &\!   0      &\!    0      &\!   -0.438 &\!    0      &\!    0      &\!   -0.074  &\!    0      &\!    0      \\
    0      &\!   0      &\!    0      &\!    0     &\!    0      &\!   -0.074  &\!    0      &\!    0      &\!    0.188  \\
    0      &\!   0      &\!    0      &\!    0     &\!   -0.073  &\!    0      &\!    0      &\!    0      &\!    0.109  \\
    0.645 &\!    0      &\!    0      &\!    0     &\!    0      &\!    0      &\!    0.188  &\!    0.109  &\!    0      \\
\end{array}
\normalsize
\right).
\end{equation*}
\normalsize
This indicates that a perturbation on the weights $a_{19}, a_{12}, a_{34}, a_{45}$ and $a_{46}$ has the highest
influence on the gap $\vp_1 - \vp_4$.

\subsection*{Modifying the full pattern of $A$} \label{sec:subf}

The choice $\eps=0.029$ with the natural initial value for the ODE \eqref{Pact-ode} given by $E(0)=-G_0$ turns out
to give interesting results.
After integration of the ODE \eqref{Pact-ode}, we get the  
Perron eigenvector of the (optimally) perturbed matrix $A + \eps E(\eps)$,
\begin{equation*}
\widetilde \vp= \small \left( 
\begin{array}{ccccccccc}
    \mathbf{0.4742} &
    0.2695          &
    0.2722          &
    \mathbf{0.4753} & 
    0.2284          &
    0.2551          &
    0.2848          &
    0.2044          &
    0.4029         \\
\end{array}
\right)^\top
\normalsize
\end{equation*} 
which shows an inversion of the ranking of the leading nodes, in fact
\begin{equation} \label{eq:exgap0}
\widetilde \vp_1 - \widetilde \vp_4 = -0.0011\ldots \qquad \Longrightarrow \qquad
\varrho_2 (\G) < 0.029,
\end{equation}
where we recall that
$
\varrho_2 (\G) =\: \min_{\widetilde{\G}} {\rm dist} \left( \mathcal G, \widetilde{\mathcal{G}} \right) = \min_{\widetilde A} \| A - \widetilde A\|,  
$
under the condition that the $2$ largest entries of the Perron eigenvector of $\widetilde A$ coalesce 
and under the constraints $\widetilde a_{ij} \ge \delta$ for $(i,j)\in\mathcal{E}$
and  $\widetilde a_{ij}=0$ for $(i,j)\notin \mathcal{E}$.
It is also interesting to observe that the gap between $\widetilde \vp_4$ and $\widetilde \vp_9$ on the other hand
has increased.

A further tuning of the parameter $\eps$ in the outer iteration provides the value $\eps^* \approx 0.0279064$ as the one leading to
coalescence (from a numerical point of view, coalescence  was declared as soon as the condition
$|\widetilde \vp_1 - \widetilde \vp_4| \le {\rm tol}$ with ${\rm tol} = 10^{-5}$ was satisfied). We obtain in fact
\small
\begin{equation*}
\widetilde \vp \approx \left( 
\begin{array}{ccccccccc}
   \mathbf{0.4745} &
   0.2696 &
   0.2717 &
   \mathbf{0.4745} &
   0.2279 &
   0.2546 &
   0.2851 &
   0.2045 &
   0.4038 
\end{array}
\right)^\top.
\end{equation*} 	
\normalsize

	

Based on an heuristic argument, one may consider a perturbation $E$ whose only nonzero entries are
$e_{12} = e_{21} = e_{14} = e_{41} = e_{19} = e_{91} = -1/\sqrt{6}$ (which guarantess the unit Frobenius norm of $E$), aiming to 
decrease the weights on the edges which connect the vertex $v_1$, which has highest centrality score, 
to the rest of the graph. 
With a perturbation again of magnitude $\eps=0.029$, this would give the Perron eigenvector of $A + \eps E$, 
\begin{equation*}
\overline \vp \approx \small \left( 
\begin{array}{ccccccccc}
    \mathbf{0.4689} &
    0.2640          &
    0.2624          &
    \mathbf{0.4543} &
    0.2229          &
    0.2523          &
    0.3047          &
    0.2161          &
    0.4270          
\end{array}
\right)^\top
\normalsize
\end{equation*} 		
which preserves the original unperturbed ranking and is characterized by a gap of $0.0146$, 
significantly different (and with the opposite sign) from the one in \eqref{eq:exgap0}.


\begin{figure} 
\begin{center}
\begin{tikzpicture}[scale=1.7, auto,swap]
    \foreach \pos/\name in {{(0,2)/v_1}, {(1,-1)/v_3}, {(4,1)/v_6}, {(6,2)/v_9},
                            {(0,0)/v_2}, {(4,-1)/v_5}, {(2.5,0)/v_4}, {(5,0)/v_7}, {(6,-1)/v_8}}
        \node[vertex] (\name) at \pos {$\name$};
    \foreach \source/ \dest /\weight in {v_2/v_1/0.160, v_4/v_1/0.270, v_9/v_1/0.240, 
                                         v_3/v_2/0.280, 
																				 v_4/v_3/0.150,
																				 v_5/v_4/0.180,v_6/v_4/0.180,
                                         v_8/v_5/0.180,
																				 v_7/v_6/0.180,
                                         v_9/v_7/0.280,
																				 v_9/v_8/0.180,
																				 v_2/v_3/\mathbf{0.277},
																				 v_5/v_8/\mathbf{0.220},
																				 v_6/v_7/\mathbf{0.216},
																				 v_7/v_9/\mathbf{0.209},
																				 v_8/v_9/\mathbf{0.137}}
        \path[edge] (\source) -- node[weight] {$\weight$} (\dest);
        \foreach \vertex / \fr in {v_4/2} 
				\node[selected vertex] at (\vertex) {$\vertex$};
				\foreach \vertex / \fr in {v_1/1} 
				\node[second selected vertex] at (\vertex) {$\vertex$};
\end{tikzpicture}
\end{center}
\caption{Modified graph of Figure \ref{fig:2} acting only on a subset of the edges 
(modified weights on each edge appear in bold) \label{fig:p}}
\end{figure}
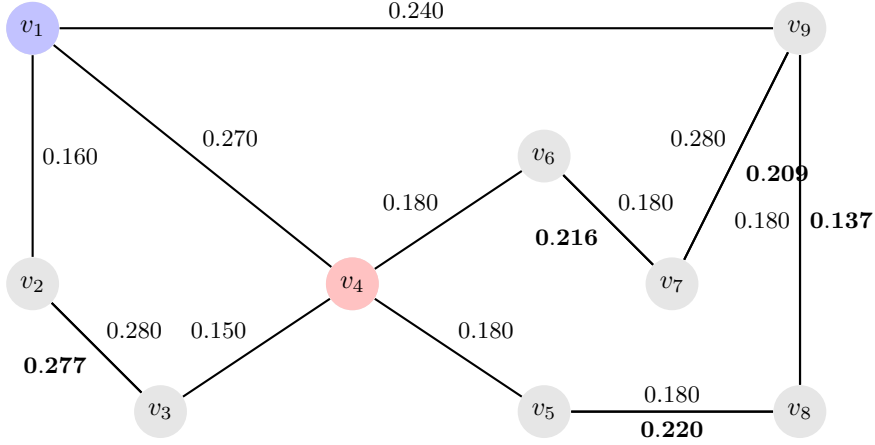

\subsection*{Modifying only a subpattern of $A$} \label{sec:subp}

As already mentioned in the discussion preceding Definition \ref{def:adm},  it is possible to consider, instead of the whole set of edges, a subset ${\mathcal{E}}' \subset \mathcal{E}$
and proceed exactly as discussed by simply replacing the projection $P_{\mathcal{E}}$ with $P_{\mathcal{E}'}$.
As an example we consider again the simple undirected graph with $9$ vertices of Figure 
\ref{fig:2} and modify only the edges not incident to the two nodes $v_1$ and $v_4$.
The results so obtained are illustrated in Figure \ref{fig:p}.
We compute the Perron eigenvector $\vp$  
\small
\begin{equation*}
\vp \approx \left( 
\begin{array}{ccccccccc}    
   \mathbf{0.4913} &
   0.2916 & 
   0.2856 &
   \mathbf{0.4913} &
   0.2406 &
   0.2608 &
   0.2433 &
   0.1894 &
   0.3601 
\end{array}
\right)^\top.
\end{equation*}
\normalsize
and the distance $\eps^* = 0.1407018$,

\subsection{Numerical integration of the gradient system}
\label{sec:num}

We describe here the $k$-th step of numerical integration of \eqref{Pact-ode}.
Starting from $E_k \approx E(t_k)$ feasible matrix of unit norm, we apply a step
of length $h$ of the explicit Euler method,
\begin{equation} \label{eq:update}
\tilde E_{k+1} = E_k + h \Pact \Bigl( {}-G_{\eps}(E_k)  + \gamma_k E_k \Bigr) \quad\hbox{ with }\quad
\gamma_k = \Big\langle  \Pact E_k, \Pact G_{\eps}(E_k) \Big\rangle 
\end{equation}
where $G_\eps$ has been computed according to \eqref{eq:grad}.

We are not interested in following the exact solution of \eqref{Pact-ode} but mainly
to decrease the functional $F_\eps$. 
This is equivalent to a gradient descent 
algorithm with line-search. 
According to this, for the stepsize selection we require that 
\begin{equation} \label{eq:cond1}
F_\eps (\tilde E_{k+1}) < F_\eps (E_k).
\end{equation}
Moreover, we apply an Armijo-like rule.
For Equation \eqref{Pact-ode} we have (recalling that $\Pact$ is just the identity unless inequality constraints activate) 
\[
\frac{d}{dt} F_\eps(E(t)) \Big|_{t=t_k} 
= -\eps \Bigl( \| \Pact G_\eps(E_k) \|^2 - \langle \Pact G_\eps(E_k), \Pact E_k \rangle ^2 \Bigr) := - g_k \le 0.
\]
We accept the result of the step with step size $h$ and set $E_{k+1} = \tilde E_{k+1}$ 
if \eqref{eq:cond1} holds, otherwise we reduce the stepsize to $h/\theta$ for some fixed $\theta>1$,
and repeat it.
Moreover, if 
\begin{equation} \label{eq:cond2}
F_\eps (\tilde E_{k+1}) \ge  F_\eps (E_k) - (h/\theta) g_k,
\end{equation}
then we reduce the step size for the next step to $h/\theta$. 
If the step size has not been reduced in the previous step, we then increase the 
step-length to be used in $[t_{k+1},t_{k+2}]$.

In a numerical solution of the gradient system (for example using the Euler-based scheme \eqref{eq:update}), 
we have to monitor the sets of edges where $a_{ij}+\eps e_{ij}=\delta$ and  among those we must further 
track those edges where the sign of $-g_{ij}+\gamma e_{ij}$ changes (where $g_{ij}$ is the $(i,j)$-th 
entry of the matrix $G_{\eps}(E)$ and $e_{ij}$ is the $(i,j)$-th entry of the matrix $E$).
When the active set $\{ (i,j)\,:\, a_{ij}+\eps e_{ij}=\delta \}$ 
is changed, then $\gamma$ also changes. 

The computational cost of a single step requires the computation of the Perron eigenvector
and the solution of the linear systems involving either the group inverse or,  in the undirected case,
the Moore-Penrose pseudoinverse (see \eqref{a-ls}).
In the undirected case we have to compute precisely $m$ columns of the Moore-Penrose
pseudoinverse.

\subsection{Outer iteration}
\label{sec:outer}


Due to the fact that in our experiments inequality constraints do not activate,
we restrict our analysis of the outer iteration to this setting (where $\Pact$
is the identical operator). A more elaborated explicit formula can be also 
deduced in the general case.

\subsubsection*{The case of inactive inequality constraints}

\begin{figure}[ht]
\centering
\includegraphics[width=.7\textwidth]{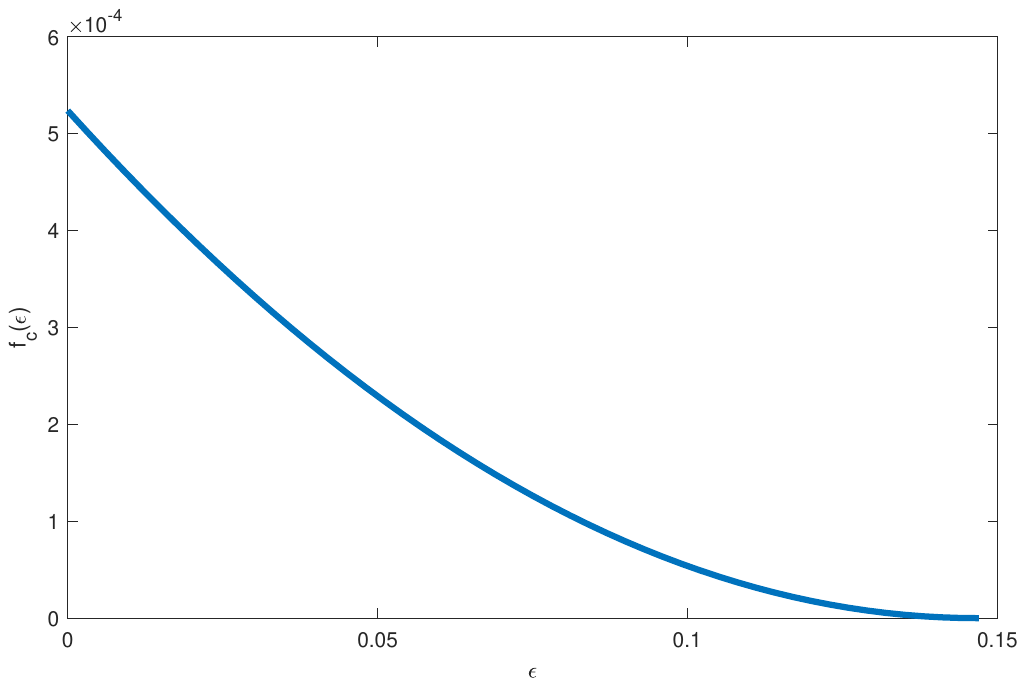}
\vspace{-3.3cm}  
\caption{Behavior of $f(\eps)$ for Example in Section \ref{sec:ill}.}
\label{fig:feps}
\end{figure}

Let $E(\eps)$ denote the minimizer of the functional $F_{\eps}$.
Generically we expect that for a given perturbation size $\eps<\epstildem$,  
$ E(\eps)$ is smooth, which means that the spectral radius of $A + \eps E$
is a simple eigenvalue. 
%
If so, then $f(\eps):=F_{\eps}(E(\eps))$ is a smooth function of $\eps$ and we 
can exploit its regularity to obtain a fast iterative method to converge to
$\epstildem$ from the left.
Otherwise we can use a bisection technique  to approach $\epstildem$.
An example is provided in Figure \ref{fig:feps}. 

The following result provides an inexpensive formula for the computation
of the derivative of $f(\eps)=F_{\eps}(E(\eps))$, which will be useful in
the construction of the outer iteration of the method. 
\begin{assumption}
We assume that $E(\eps)$ is a smooth function of $\eps$ in some interval. 
\label{assumpt}
\end{assumption}
\vspace{-0.12in}
We then have the following result.
\begin{lemma}
\label{lem:der}
Under Assumption~{\rm \ref{assumpt}}, the function $f(\eps)=F_{\eps}(E(\eps))$ 
is differentiable and its derivative equals (with ${\phantom{a}'}= d/d\eps$)
\begin{equation}
f'(\eps) =  - \|  G_{\eps}(E(\eps)) \|  
\label{eq:derFdeps}
\end{equation}
\end{lemma}
\begin{proof} 
The  conservation of $\|  E(\eps) \|=1 $ implies
\[
\Big\langle E(\eps), E'(\eps)\Big\rangle = 0.
\]
Differentiating $f(\eps)=F_{\eps}\bigl( E(\eps) \bigr)$ with respect to $\eps$ we obtain
\begin{eqnarray}
\hskip -9mm
f'(\eps) 
=\Bigl\langle G_{\eps}(E(\eps)), E(\eps)+\eps E'(\eps) \Bigr\rangle 
=\Bigl\langle G_{\eps}(E(\eps)), E(\eps) \Bigr\rangle.
\label{eq:dFdeps}
\end{eqnarray}
By means of the property no.~3. of minimizers, stated in Theorem~\ref{thm:stat}, we can conclude that
$E(\eps) \propto G_{\eps}(E(\eps))$, which yields formula \eqref{eq:derFdeps}.
\end{proof}

%


\subsubsection*{Newton method}
For $\eps < \epstildem$, we make use of the standard Newton iteration
\begin{equation}
\eps_1 = \eps_o - \frac{f(\eps_o)}{f'(\eps_o)},
\label{eq:Newton}
\end{equation}
to get an updated value $\eps_1$ from $\eps_o$.

In a numerical code it is convenient to couple the Newton iteration \eqref{eq:Newton} 
with the bisection method, to increase its robustness.
We make use of a tolerance which allows us to distinguish numerically whether $\eps < \epstildem$, 
in which case we may use the derivative formula and perform the Newton step, or $\eps > \epstildem$, 
in which case we just bisect the current interval.
\vspace{-3.3cm}  

%
%



\section{Numerical tests}
\label{sec:numexp} 


In this section we present numerical experiments illustrating the behavior of our algorithm on a
few undirected and directed sparse weighted graphs arising in applications.   We also perform a
comparison with more standard approaches for solving constrained optimization problems.

All tests are made on a Surface Book 3 computer, under the {\sc Matlab} version 23.2.0.2365128 (R2023b).

\subsection{Example 1, matrix $A_{27}$}

\begin{figure}[ht]
\centering
\includegraphics[width=.48\textwidth]{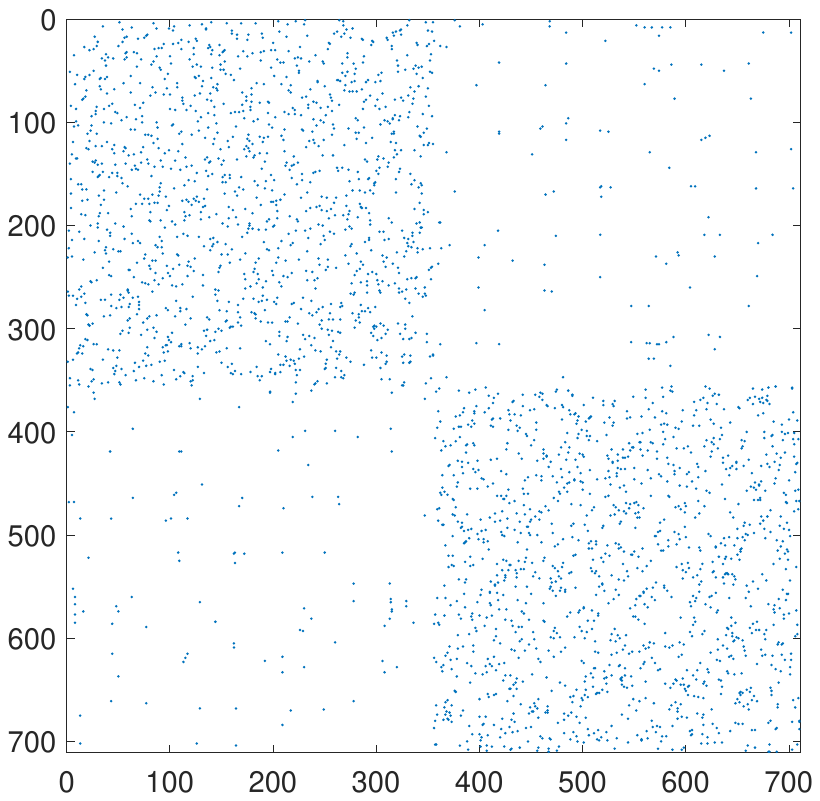} 
\includegraphics[width=.48\textwidth]{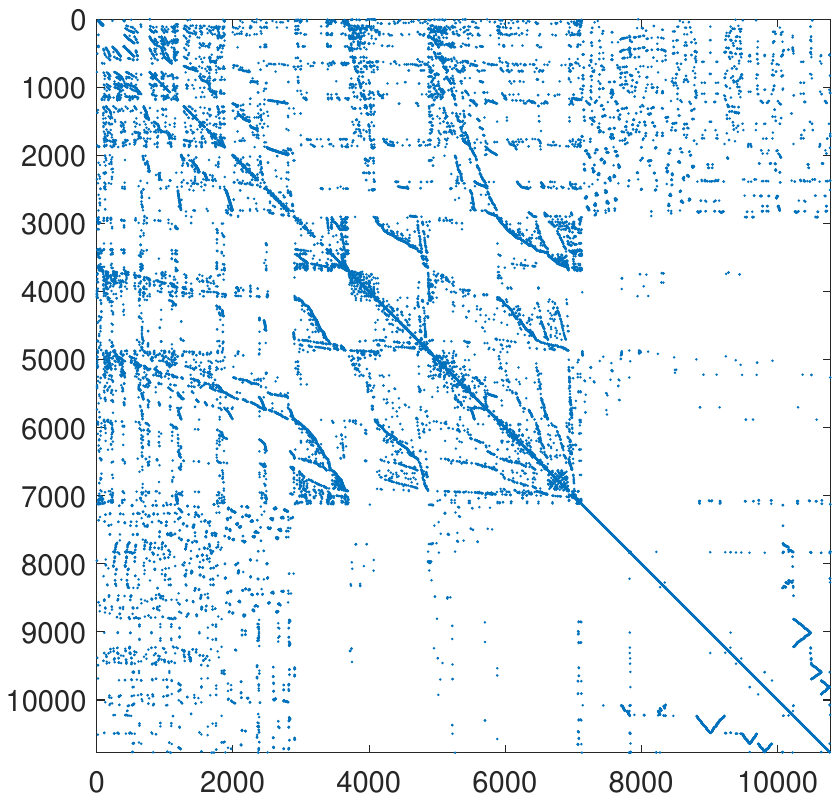} 
\vspace{-2cm}
\caption{Patterns of the adjacency matrix for the graph example $A_{27}$ (left) and for the graph example \emph{nopoly} from Gaertner collection (right)}
\label{fig:ex1_1}
\label{fig:ex2_1}
\end{figure}

The considered graph is one of a sequence of graphs arising from molecular dynamics simulations aimed at 
analyzing the structure of water networks under different temperature and pressure conditions (see \cite{Faccio_et_al}).
It consists of $710$ nodes (corresponding to oxygen atoms), most of which ($613$) have degree $4$. Its pattern is shown in Figure \ref{fig:ex1_1}.
\begin{figure}[ht]
\centering
\includegraphics[width=.48\textwidth]{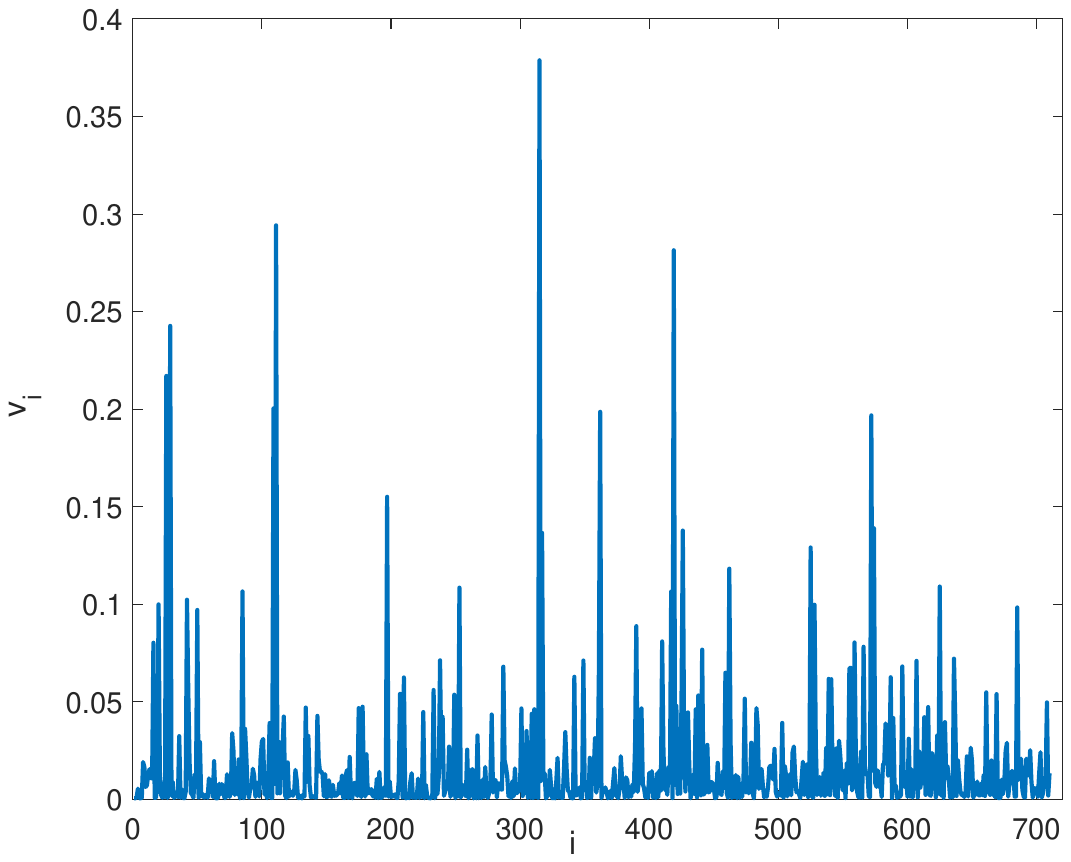} \includegraphics[width=.48\textwidth]{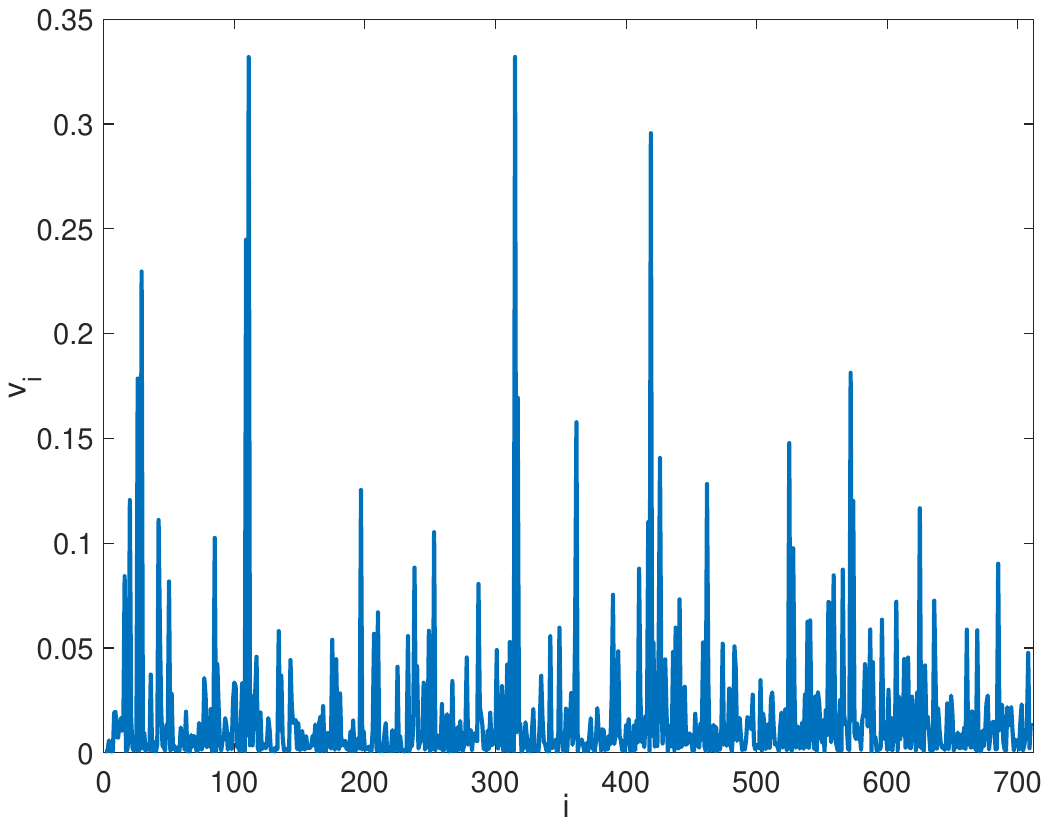}
\vspace{-2cm}
\caption{Distribution of the entries of Perron eigenvector for the example $A_{27}$. Left picture: original matrix. Right picture: perturbed matrix
(for $m=2$).}
\label{fig:ex1_2}
\end{figure}
In Figure \ref{fig:ex1_2} we show the entries of the Perron eigenvector of the original (left picture) and perturbed (right picture),
obtained for $m=2$.

The five leading entries of the Perron eigenvector of $A$ are the following:
\begin{eqnarray*}
&&
\vp_{26}  = 0.21705332, \quad
\vp_{29}  = 0.24269500, \quad 
\vp_{419} = 0.28145420, 
\\
&&
\vp_{111} = 0.29422909, \quad 
\vp_{315} = 0.37879622.
\end{eqnarray*}
After scaling the matrix $A$ to have Frobenius norm $1$, we apply the algorithm.

The norm of the perturbation which determines the coalescence of the two main entries of the Perron eigenvector is
$\widetilde \eps_2 = 0.0092654$, which is below $1\%$. The two leading entries after the perturbation are
\[
\vp_{111}=0.{\bf 331}01615, \qquad 
\vp_{315}=0.{\bf 331}36357.
\]
The number of outer iterations is $16$. The CPU time is about 1.23 seconds. 
\begin{table}[hbt]
\begin{center}
\begin{tabular}{|l|l|l|l|}\hline
  $k$ & $\eps_k$ & $f_c(\eps_k)$ & $\#$ eigs \\
 \hline
\rule{0pt}{9pt}
\!\!\!\! 
  $0$         & $0.005$ &  $0.000375176389506$ & $25$  \\
	$1$         & $0.005136626814311$ & $0.000351719368954$ & $2$  \\
	$2$         & $0.008469968038737$ & $0.000013522265277$ & $34$  \\
	$3$         & $0.008489282220273$ & $0.000012884575588$ & $2$  \\
	$6$         & $0.009184673253577$ & $0.000000185817345$ & $2$  \\
	$10$        & $0.009212325217551$ & $0.000000093285281$ & $1$  \\ 
 \hline
\end{tabular}
\vspace{2mm}
\caption{Computation of  $\eps_k$, $f_c(\eps_k) = F_\eps^c(E_k)$ 
and number of eigenvalue computations in the inner iterations for Example 1, matrix $A_{27}$
with $m=2$.  
 \label{tab:ex1}}
\end{center}
\end{table}

Note that the relative variations of the two entries corresponding to the nodes with highest centrality are
$\delta_{111} = +0.1281$ and $\delta_{315} = -0.1235$.

We also consider coalescence of the first $m=3$ entries, for which we obtain a perturbation of norm $\widetilde \eps_3 = 0.010718306874804$
with the $3$ largest entries of Perron eigenvector $\vp$ given by  
\small
\[
 \vp_{419} =  0.{\bf 3000}0085, \quad
 \vp_{315} =  0.{\bf 3000}2148, \quad
 \vp_{111} =  0.{\bf 3000}2174.
\]
\normalsize
To obtain $4$ digits of precision we required a smaller tolerance to step the outer iteration, which determined a CPU time
of about $9.33$ seconds.
Finally we set $m = 4$ for which we obtain 
$\widetilde \eps_4 = 0.014389966738210$, with the largest entries of the Perron eigenvector $\vp$ being
\small
\[
\vp_{109}  = 0.{\bf 2317}9527, \quad
\vp_{419} =  0.{\bf 2318}0413, \quad
\vp_{111} =  0.{\bf 2318}0684, \quad
\vp_{315} =  0.{\bf 2318}0693 .
\]
\normalsize
The CPU time is about $16$ seconds.
It is interesting to note the the node $109$ was not at the highest rank but replaced the node $29$ after the perturbation.
\smallskip

\subsubsection*{Comparison with {\sc Matlab} routines} 

In a first experiment 
we have applied the {\sc Matlab} routine {\tt fmincon} to the problem formulation \eqref{eq:Omega}--\eqref{eq:nearness}
for the case $m=2$. Here we use {\tt fmincon} as a black box, without providing the analytic form of the gradient (the
derivation of which is among the results of this paper) to the routine.  
Using standard default values we have not been able to obtain a correct solution, hence it has been necessary to 
increase the number of function evaluations allowed. Table \ref{tab:Mat1} shows the results of the routine with an interior point
algorithm, which is the one that performs best.
The resulting solution is very close to the one we found, with coalescent entries $\# 111$ and $\# 315$,
\[
\vp_{111}=0.{\bf 33}2\ldots, \qquad 
\vp_{315}=0.{\bf 33}2\ldots.
\]

In Table \ref{tab:Mat1}, $k$ stands for the iteration number, $x$ is the vector containing the nonzero entries of the
upper triangular part of the weighted adiacency matrix $A$,  $g(x)$ is the constraint, $\vp_1 - \vp_2$ being $\vp(x)$ the
Perron eigenvector of the matrix $B = A - X$, where $X$ is the perturbation associated to $x$. 
Finally, $f(x) = \sqrt{2} \| x \| = \| X \|_F$ and $\# fe$ indicates the number of function evaluations. 
The overall CPU time is about $80$ minutes.
\begin{table}[hbt]
\begin{center}
\begin{tabular}{|l|l|l|l|}\hline
  $k$ & $\#$ fe & $f(x_k)$ & $g(x_k)$ \\
 \hline
\rule{0pt}{9pt}
\!\!\!\! 
0  & 1.4750e+03  & 6.9396e-04 & 8.4470e-02 \\ 
1  & 2.9500e+03 & 2.4405e+01 & 7.6330e-02 \\
6  & 1.0325e+04 & 3.8498e+01 & 4.3180e-04 \\
11 & 1.7700e+04 & 7.3944e+00 & 2.0100e-03 \\
21 & 3.2450e+04 & 6.7381e-01 & 2.5690e-03 \\
31 & 4.7206e+04 & 6.6791e-01 & 1.7290e-06 \\
41 & 6.1970e+04 & 5.1159e-02 & 2.3090e-02 \\
51 & 7.6743e+04 & 4.9650e-02 & 2.9350e-04 \\
61 & 9.1510e+04 & 4.0801e-02 & 1.0950e-02 \\
71 & 1.0629e+05 & 1.4522e-02 & 4.7660e-03 \\
81 & 1.2109e+05 & 9.5850e-03 & 3.9870e-05 \\
91 & 1.3589e+05 & 9.2726e-03 & 3.5810e-10 \\
96 & 1.4329e+05 & 9.2726e-03 & 5.8570e-10 \\
\hline
\end{tabular}
\vspace{2mm}
\caption{Reported behavior of the {\sc Matlab} algorithm.  
 \label{tab:Mat1}}
\end{center}
\end{table}
Observing Table \ref{tab:Mat1}, we see that in the beginning $x$ gets quite large with respect to its final value
($3.8 \cdot 10^1$ versus $9.27 \cdot 10^{-3}$). This very likely implies activation of inequality constraints and
makes the convergence process slower. For this reason, trust region techniques would likely improve the process.
The initial value $x_0$ provided to the algorithm is quite small (as $f(x_0)$ indicates), but the algorithm makes 
a large jump to a vector $x_1$ of much higher norm, which is due to the violation of the constraint. Starting from a
consistent value would help, but requires additional computation.
\medskip

Next we describe the results of a second experiment, with analytic gradients provided to the algorithm.
Following the notation in Section \ref{sec:consprob}, we have provided to the algorithm, the exact formula 
for the gradient of the constraint $u_{i_1} - u_{i_2} = 0$, which is obtained by the same arguments used in Lemma \ref{lem:a-dot} 
and is proportional to 
\[
P_{\mathcal{E}} ( \sym(B^\dagger (u_{i_1} - u_{i_2}) \vp^\top ).
\]
The norm of the perturbation which determines the coalescence of the two main entries of the Perron eigenvector is
$\widetilde \eps_2 = 0.01023234$, which is about $1\%$. The two leading entries after the perturbation are
\[
\vp_{315} \approx \vp_{636}= 0.20245843.
\]
This indicates a different local optimum with respect to the previous experiment (which was the same as the one we found).
However, they are very close and both indicate a strong sensitivity of the ranking for this graph.
Table \ref{tab:Mat2} shows the results of the run.
\begin{table}[hbt]
\begin{center}
\begin{tabular}{|l|l|l|l|}\hline
  $k$ & $\#$ fe & $f(x_k)$ & $g(x_k)$ \\
 \hline
\rule{0pt}{9pt}
\!\!\!\! 
    0  &  1.475e+03   &  1.002577e-03  &  8.433e-02  \\
    1  &  2.950e+03   &  2.329900e+01  &  5.114e-02  \\
   21  & 3.2464e+04   &  4.892883e+00  &  5.230e-03  \\
   41  & 6.1970e+04   &  2.175929e-01  &  5.546e-03  \\
   61  & 9.1493e+04   &  1.814369e-01  &  6.202e-06  \\
   81  & 1.21020e+05  &  1.786729e-02  &  1.009e-03  \\
  101  & 1.50620e+05  &  1.025323e-02  &  9.777e-08  \\
  121  & 1.80195e+05  &  1.023304e-02  &  3.949e-07  \\
  141  & 2.09771e+05  &  1.023236e-02  &  1.121e-07  \\
  147  & 2.18661e+05  &  1.023233e-02  &  2.660e-09  \\
\hline
\end{tabular}
\vspace{2mm}
\caption{Reported behavior of the {\sc Matlab} algorithm provided by the analytic gradient of the
constraint. 
 \label{tab:Mat2}}
\end{center}
\end{table}
The overall CPU time is about $351$ seconds, which shows a significant improvement with respect to 
the previous experiment, although still far from the CPU time required by our approach.

\subsection{Example $2$, matrix {\em nopoly}}

This matrix has been contributed to the TAMU Sparse Matrix Collection by Dr.~Klaus Gaertner, Weierstrass Institute, Berlin\footnote{See https://sparse.tamu.edu/Gaertner}.
The matrix, described as the adjacency matrix of a weighted undirected graph, is symmetric. However, some of the weights are negative. 
To make it non-negative we consider its absolute value. 
The corresponding graph
consists of $n=10774$ nodes. 
Its pattern is shown in Figure \ref{fig:ex2_1}.  We do not perform a comparison with {\sc Matlab}'s optimization
tools since the size of the matrix makes this prohibitively expensive. 

The five leading entries of the Perron eigenvector of $A$ are the following:        
\begin{eqnarray*}
&& \vp_{1155} =   0.12197226, \quad 
\vp_{3391} =   0.12460486, \quad
\vp_{3002} =   0.12460960, 
\\
&& \vp_{3007} =   0.12502517, \quad
\vp_{3397} =   0.12559705 .
\end{eqnarray*}
After normalizing the matrix $A$ in the Frobenius norm, we apply the algorithm with the aim of obtaining the largest $m$
entries of the Perron eigenvector to match to a $3$-digit accuracy.

First we apply the algorithm with $m=2$.
The method starts with $\eps_0 = 10^{-3}$ and sets the tolerance for the stopping criterion to $10^{-6}$. 
It makes $3$ outer iterates and the final value of the functional is $8 \cdot 10^{-12}$.
The norm of the perturbation which determines the matching of the $2$ leading entries of the Perron eigenvector is
$\widetilde \eps_2 = 0.00383633$, which is about $0.00008\%$. The computed $2$ leading entries are the following:   
\[
\vp_{3007} =   0.{\bf 12558}011, \quad
\vp_{3397} =   0.{\bf 12558}553.
\]
The CPU time is about $19.2$ seconds. 

Then we apply the algorithm with $m=3$.
The method starts with $\eps_0 = 10^{-3}$ and sets the tolerance for the stopping criterion to $10^{-10}$. The method makes $8$ outer iterates.
The norm of the perturbation which determines the coalescence of the $3$ leading entries of the Perron eigenvector is
$\widetilde \eps_3 = 0.003909011$, which is about $0.00008\%$. The computed $3$ leading entries are the following.   
\[
\vp_{3007} =   0.{\bf 12550}017, \quad
\vp_{3002} =   0.{\bf 12550}316, \quad
\vp_{3397} =   0.{\bf 12550}872 
\]
which coincide to a $5$-digit approximation.
The CPU time is about $30.8$ seconds. 

Next we set $m=5$.
The norm of the perturbation which determines the coalescence of the $5$ leading entries of the Perron eigenvector is
$\widetilde \eps_5 = 0.0161$, which is about $0.00016\%$. 
The new leading entries follow:  
\begin{eqnarray*}
&& \vp_{1155} =   0.{\bf 124}48751, \quad 
\vp_{3391} =   0.{\bf 124}51103, \quad
\vp_{3007} =   0.{\bf 124}53422, 
\\
&& \vp_{3002} =   0.{\bf 124}60202, \quad
\vp_{3397} =   0.{\bf 124}63431 .
\end{eqnarray*}
The CPU time is about $321.8$ seconds.

\begin{figure}[ht]
\centering
\includegraphics[width=.6\textwidth]{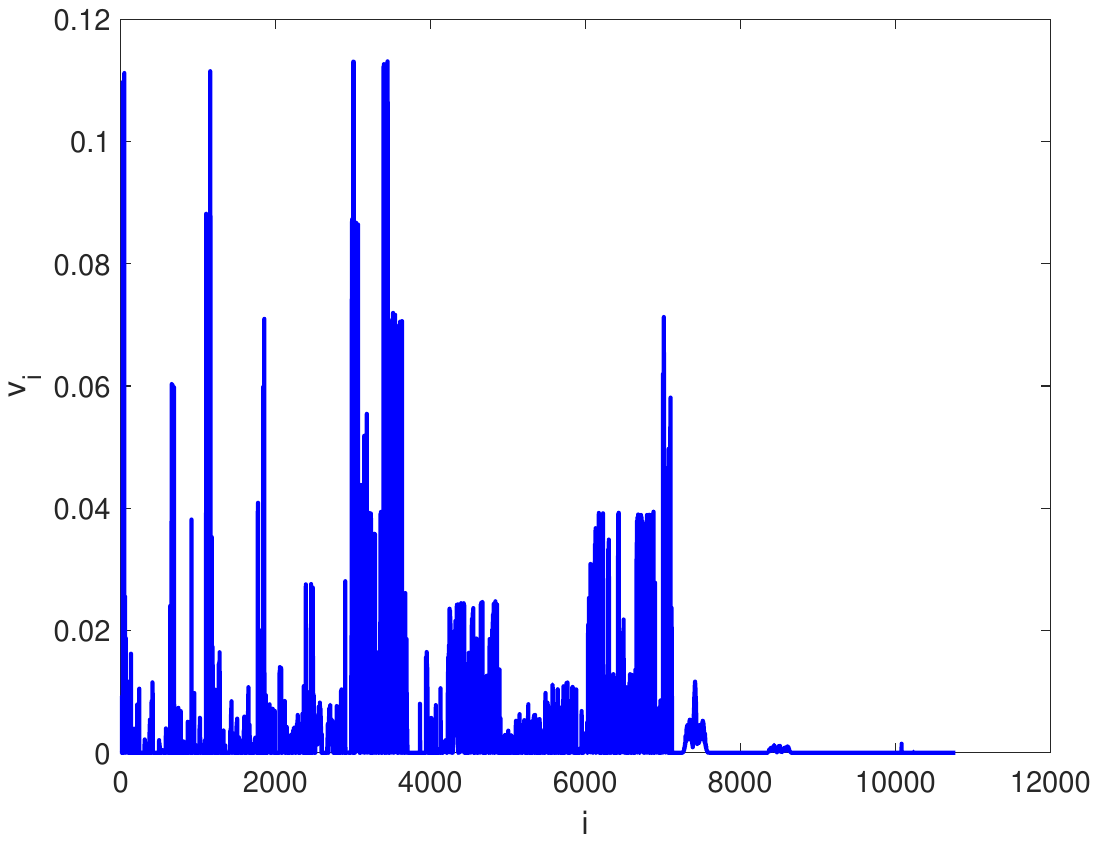} 
\vspace{-2cm}   
\caption{Distribution of the entries of Perron eigenvector of the matrix obtained for $m=10$ for the example \emph{nopoly}  from Gaertner collection.}
\label{fig:ex2_10}
\end{figure}

Finally, in Figure \ref{fig:ex2_10} we show the entries of the Perron eigenvector of the  perturbed matrix obtained setting $m=10$.
\vspace{-1cm}
\newpage

\section{Conclusions} \label{sec:concl}
In this paper we have presented a methodology for finding a minimal perturbation
to the adjacency matrix of a weighted (sparse) graph that causes the leading entries of the Perron
eigenvector of the perturbed matrix to coalesce. The perturbation is constrained so as to preserve
the sparsity of the original matrix. One of the motivations for this work was the question of the reliability
and robustness of eigenvector centrality: if a small perturbation causes coalescence of the
leading entries of the Perron eigenvector, this centrality measure cannot be considered reliable in the
presence of uncertainties or noise, especially when the leading entries of the perturbed eigenvector 
differ from the original ones  and the size of the perturbation is within the uncertainty limits.   
  Our methodology easily allows to incorporate additional constraints; for example,
only a subset of the graph's edge weights may be allowed to change. This feature can be used to identify which
edges play an important role in influencing the ranking of nodes in the graph.
The proposed algorithm appears to perform  favorably when compared with standard constrained optimization methods
 on graphs arising from  realistic applications, and is fairly efficient in spite of its complexity.  

%
%

\bigskip
\bigskip

\end{document}